\documentclass[12pr0]{article}
\usepackage{amssymb}
\usepackage{amsmath}
\usepackage{boxedminipage}
\usepackage{graphicx}
\usepackage{float}
\usepackage[nodayofweek,level]{datetime}
\usepackage{geometry}
\begin{document}
    \title{Characteristic-Sweeping for Quick Solution Generation and Shock Prediction of 1D Scalar Conservation Law}
    \author{Gaojin He}
    \date{August 8, 2018}
    \maketitle
    \centerline{Department of Applied Mathematics and Statistics}
    \centerline{Stony Brook Univeristy}
\bigskip\bigskip\bigskip
\begin{abstract}
\noindent As we known for simple first order partial differential equation $u_t+(G(u))_x=0$ (In this paper, function $G$ is smooth enough), there are several popular methods to be applied, such as Godunov Scheme with limiters, WENO and DG. Notice that the main problem is these method all require iteration of time step for all $x$ span. Because of this large scale iteration, often one spends decent time on running them, and uses large memory space to store the dynamic solution. \\

\noindent This paper is to establish a new possible numerical method to quickly solve first order one-dimensional scalar hyperbolic equation. This method bases on the very original fact of the property of first order PDE--characteristics, and significantly reduces the reliance of time step iteration for solving this PDE numerically. At the same time, this method can be coded into algorithm as a super memory-saving precondition to store the solution from $t=0$ to $t=T$, or any time span. Moreover, it might be possible to generate this method to solve the PDE in higher dimension case. However of course, one should first understand how this works in 1D.\\

\noindent If we allow muti-value solution, then this equation can be easily solve completely through characteristic. The result can be written as $F(x,t,u)=0$ mathematically and forms a surface in $\mathbb{R}^3$. However this is not the result we want. The solution we need must be single-value, and physically satisfies entropy condition. Because of this restriction, the solution may contain shocks, rarefaction waves, or both, which are the main difficulty we encounter if we continue using characteristic to solve the equation.\\

\noindent In this paper, I will illustrate how to overcome this difficulty by establishing the equations for shocks and predict them. Then we can see how one can sweep the characteristic to get the solution at any time $t=T$. I will start first from the prerequisite condition of the weak solution, then base on this establish the equation of shock curve, and finally show how to numerically predict them and form an algorithm.\\
\end{abstract}
\newpage
\tableofcontents
\newpage

\section{Prerequisite}
This is a section for those who are not familiar with the problem we are going to discuss. If you have already known the content of this section, you can skip to the next section to see how the equation of shock wave be established.
\subsection{Characteristic on Multi-Value Parameterized Solution}
We first to see how normally we solve\\\\
$$u_t+(G(u))_x=0 \eqno{(1.1.1)}$$\\\\
with $u(x,0)=f(x)$and for $t\ge0$ by using characteristic. This is equivalent to the problem of finding a function $F$ such that $F(x,t,z)=0$ by denoting $z=u(x,y)$, and the solution is indeed the level curve on the graph of function $F$ in $\mathcal{R}^4$. Since we know the solution is of the form $z=u(x,y)$, we can let\\\\
$$F \triangleq F(x,t,z)=u(x,t)-z \eqno{(1.1.2)}$$\\\\
Define the parameterized characteristic $\mathcal{C}(s)=(x(s),t(s),z(s))$ and substitute into the solution $ F(x,t,u)=0$ we have\\\\
$$u(x(s),t(s))-z(s)=0 \eqno{(1.1.3)}$$\\\\
Take the derivative with respect to s, then we have\\\\
$$u_x\frac{\mathrm{d}x}{\mathrm{d}s}+u_t\frac{\mathrm{d}t}{\mathrm{d}s}=\frac{\mathrm{d}z}{\mathrm{d}s} \eqno{(1.1.4)}$$\\\\
Note that $u_t+(G(u))_x=0$ is equavalent to $u_t+G'(u)u_x=0$, compare with (1.1.3) we have\\
$$\frac{\mathrm{d}t}{\mathrm{d}s}=1 \qquad \frac{\mathrm{d}x}{\mathrm{d}s}=G'(z) \qquad \frac{\mathrm{d}z}{\mathrm{d}s}=0 \eqno{(1.1.5)}$$\\\\
We restrict characteristic start from initial condition with s=0 \quad $i.e.\quad \mathcal{C}(0)=(x(0),t(0),z(0))=(x^*,0,f(x^*))$ for any point $x^*$. This is a first order ODE system with exact solution\\\\
$$t=s \qquad x=G'(f(x^*))s \qquad z=f(x^*) \eqno{(1.1.6)}$$\\\\
This means the characteristic pass through initial point $(x^*,0,f(x^*))$ is\\\\
$$x=G'(f(x^*))t \qquad z=f(x^*) \eqno{(1.1.7)}$$\\\\
We can see the characteristic is a straight line passes through $(x^*,0,f(x^*))$ for any $x^*$ with no change of $z$ coordinate. This, along with (1.1.3), means along the line $x=G'(f(x^*))t$ it is always $u(x,t)\equiv f(x^*)$ . For any $x$ and $t$, the characteristic passes through $(x,t)$ intersect the plane $t=0$ at $(x-G'(z)t,0,z)$. Because of this, we must have\\\\
$$z=f(x-G'(z)t) \qquad where\quad z=u(x,t) \eqno{(1.1.8)}$$\\\\
This is the implicit multi-value solution for the problem (1.1.1). Geometrically, the set of all characteristics forms a surface which is the overall solution for any $x$ and $t$. We only care about the solution of $t\ge0$, and for any $t=T$ the multi-value solution $z=u(x,t)$ is the curve of intersection of the surface and the plane $t=T$.\newpage

\subsection{Shock for Single-Value Solution}
Recall the original problem $u_t+(G(u))_x=0$, the notation $u(x,t)$ implies we only accept single-value solution \quad $i.e$ \quad for any fix $x$ and $t$, there is only one single value of $u$. But through the characteristic we know mathematically this is not guranteed. When two characteristics intersect at time $t^*\ge0$ at $x^*$, then $u(x^*,t^*)$ has at least two different value originated from this two characteristics at $t=0$. That is, geometrically, the formation of the multi-value solution.\\\\
To avoid this, one should realize that in the case of characteristics intersect, forcing it become single-value solution means there will be a discontinuity of the solution. We call this discontinuity \textbf{Shock}.\\\\
To solve this, we have to see how this equation was derived. Intergrate (1.1.1) with respect to $x$ from a to b, we obtain\\\\
$$G(u(b,t))-G(u(a,t))+\frac{\mathrm{d}}{\mathrm{dt}}\int_a^bu(x,t)\mathrm{d}x=0 \eqno{(1.2.1)}$$\\\\
If $u(x,t)$ is continuous from $x=a$ to $x=b$, then (1.2.1) is trivial, and equivalent to (1.1.1). However, when shock forms, denote the $x$ coordinate of the shock $\xi(t)$ then we will see the following if $a<\xi(t)<b$\\\\
$$G(u(b,t))-G(u(a,t))+\frac{\mathrm{d}}{\mathrm{dt}}(\int_a^{\xi(t)} u(x,t)\mathrm{d}x+\int_{\xi(t)}^b u(x,t)\mathrm{d}x)=0 \eqno{(1.2.2)}$$\\\\
Which is
$$G(u(b,t))-G(u(a,t))+ \xi'(t)(u_l-u_r)+\int_a^{\xi(t)} u_t(x,t)\mathrm{d}x+\int_{\xi(t)}^b u_t(x,t)\mathrm{d}x=0 \eqno{(1.2.3)}$$\\\\
where\\\\
$$u_l=\lim_{x \to \xi^-(t)}u(x,t) \qquad \qquad u_r=\lim_{x \to \xi^+(t)}u(x,t )\eqno(1.2.4)$$\\\\
We let $a\uparrow\xi(t)$ and $b\downarrow\xi(t)$ then we obtain\\\\
$$\xi'(t)=\frac{G(u_r)-G(u_l)}{u_r-u_l} \eqno(1.2.5)$$\\\\
This is the \textbf{jump condition} of shock.\\\\
\par
\centerline{\includegraphics[width=18cm,height=14cm]{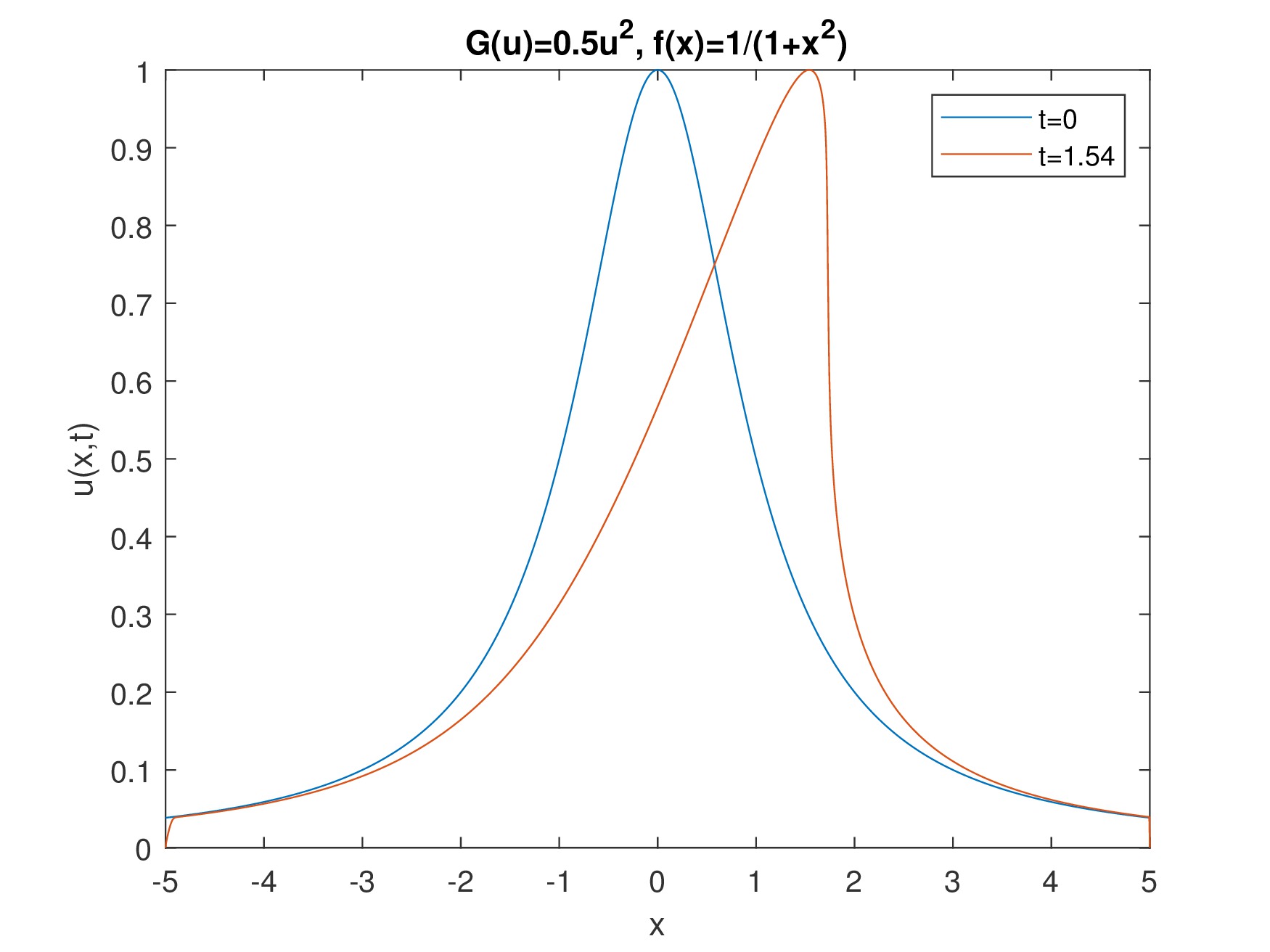}}
\centerline{Shock Forms at continuous point}
\newpage

\subsection{Rarefaction Wave}
Unlike the shock, rarefaction wave only occurs when initial condition has a jump discontinuity. In addition, if the jump discontinuity point of initial condition f(x) is $x^*$, then rarefaction wave occurs if and only if $G'(f(x^{*-}))<G'(f(x^{*+}))$. This is because in this case, for some interval $t\in(0,\delta)$ there will be no characteristic in the region\\\\
$$\{(x,t)\mid G'(f(x^{*-}))<\frac{x-x^*}{t}<G'(f(x^{*+}))\} \eqno(1.3.1)$$\\\\
In such region we must be creative and construct a solution that is single-value, and satisfies (1.1.1). One smart way is to solve the equation with treating the discontinuity point as a very short segment and take the limit. Locally, we extend the discontinuity a little bit, make it as $f(x)=u(x,0)=k(x-x^*)$. Recall (1.1.8) we have\\\\
$$z=k(x-G'(z)t-x^*) \eqno(1.3.2)$$\\\\
which is\\\\
$$G'(z)=\frac{x-x^*-\frac{z}{k}}{t} \eqno(1.3.3)$$\\\\
It is natrual to think of the jump discontinuity as the case when $k\to+\infty$. Thus, we take the limit $k\to+\infty$ at (1.3.3) and obtain\\\\
$$G'(z)=\frac{x-x^*}{t}$$\\\\
which is\\\\
$$u(x,t)=(G')^{-1}(\frac{x-x^*}{t}) \eqno(1.3.4)$$\\\\
The solution (1.3.4) is called the \textbf{rarefaction wave}. It is easy to directly verify (1.3.4) that indeed satisfies (1.1.1), and is continuous in the region (1.3.1). Moreover, (1,3,4) also guarantees the continuity of $u(x,t)$ along two rays $\frac{x-x^*}{t}=G'(f(x^{*-}))$ and $\frac{x-x^*}{t}=G'(f(x^{*+}))$.\\\\
One may notice that in the region (1.3.1), there is another way to construct the solution which is not continuous but still satisfies (1.2.5). We call this solution \textbf{rarefaction shock}.  However this mathematically correct solution physically violates the entropy condition, which means in real world this could never happen. We just  remember when it is the region (1.3.1), the rarefaction wave is the only solution we need to consider.\\\\
\par
\centerline{\includegraphics[width=18cm,height=14cm]{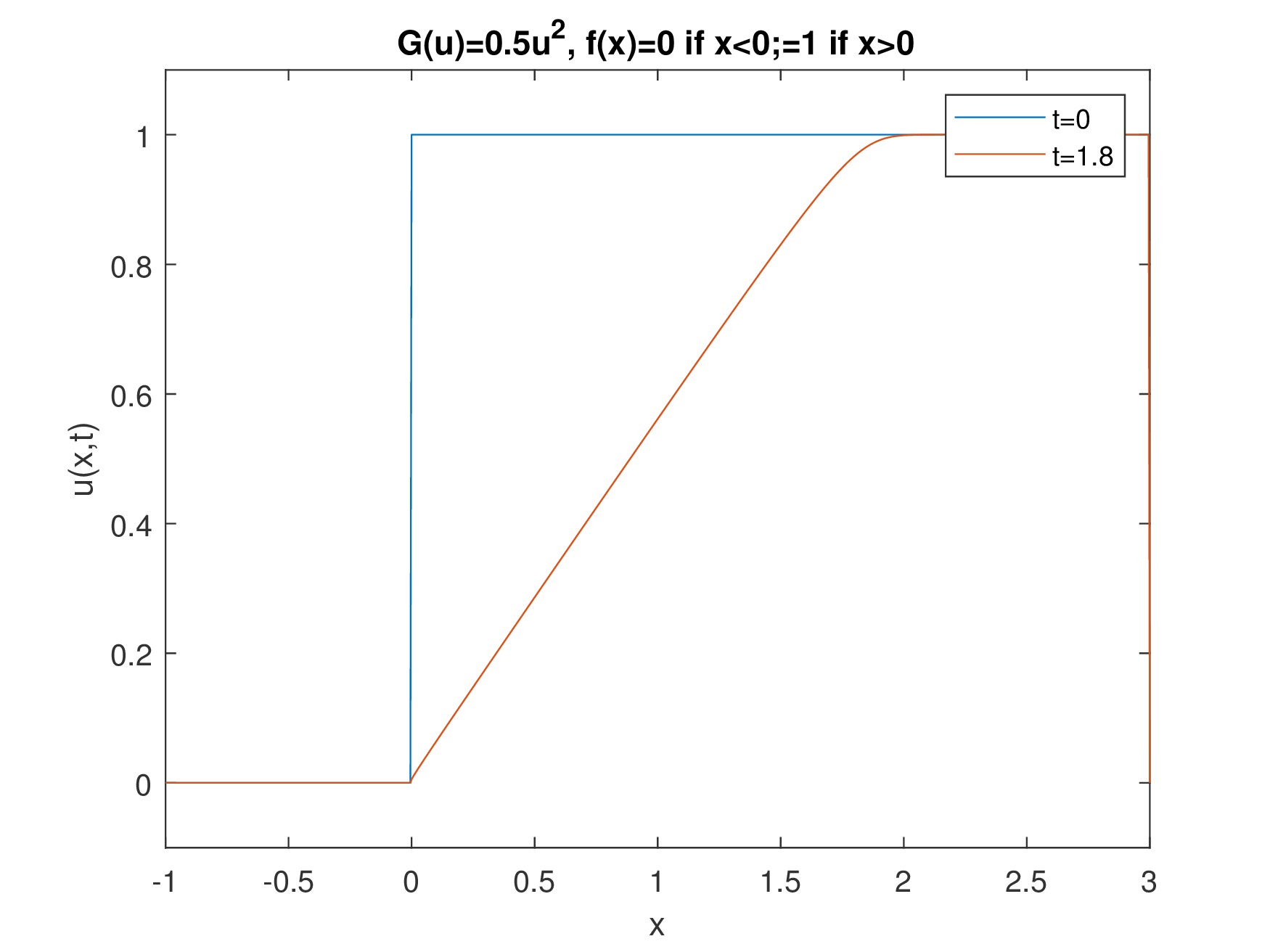}}
\centerline{Rarefaction}
\newpage

\section{Local Shock Formation}
In this section we mainly analyze those cases in which shock wave forms and, establish the equations for shock wave. One notes that shock waves may intersect and form a new shock wave. We will first deduct the shock wave equations locally, and then in the next section analyze the case when shock waves intersect.

\subsection{Local Shock Point and Break Time}
The main difficulty of using characteristic to solve the weak solution of (1.1.1) is the shock wave. Before we see how we can predict them, one must know in what case will shock form.\\\\
Naively we can write down the equations of 2 characteristic line that potentially intersect, then solve them as a system and get there intersect time $t$.\\\\
$$
\begin{cases}
G'(f(x_l))t=x-x_l\\
G'(f(x_r))t=x-x_r\\
\end{cases}
\eqno(2.1.1)$$\\\\
Which implies\\\\
$$t=-\frac{x_l-x_r}{G'(f(x_l))-G'(f(x_r))}\eqno(2.1.2)$$\\\\
There are 3 cases which yield different local smallest t from (2.1.2). We need to discuss each of them.\\\\\\
$\mathbf{Case\ 1}$\\\\
Assume $f(x)$ has a jump point $x^*$. Note that we only care about the solution of $t\ge0$, and since we have set $x_l<x_r$, by taking $x_l\uparrow x^{*-}$ and $x_r\downarrow x^{*+}$, it can be seen that if $G'(f(x^{*-}))>G'(f(x^{*+}))$ then shock will form immediately at point $x=x^*$ at time $t=0$. (Recall that if $G'(f(x^{*-}))<G'(f(x^{*+}))$ then it is rarefaction wave.)\\\\
\textbf{Definition\ 1:} \emph{A point $x^*$ is called the \textbf{first kind shock point} if $x^*$ is a jump point of $f(x)$ and $G'(f(x^{*-}))>G'(f(x^{*+}))$}.\\\\\\
$\mathbf{Case\ 2}$\\\\
Assume $f(x)$ is continuous at point $x^*$, but $x^*$ is a jump of $f'(x)$. In this case shock may occur and, it is a little bit complicated. Let us put a definition here first and we shall discuss later.\\\\
\textbf{Definition\ 2:} \emph{A point $x^*$ is called the \textbf{second kind shock point} if the following conditions are satisfied}\\\\
1 \quad \emph{$f(x)$ is continuous at $x^*$}\\\\
2 \quad \emph{$x^*$ is a jump of $f'(x)$}\\\\
3 \quad \emph{$\frac{\mathrm{d}}{\mathrm{d}x}G'(f(x^{*-}))<0$ or $\frac{\mathrm{d}}{\mathrm{d}x}G'(f(x^{*+}))<0$}\\\\
4 \quad \emph{There exist $\delta_0>0$ such that for any $\delta\in(0,\delta_0)$, at least one of the following \emph{(1) (2)} is satisfied}\\\\
\emph{\emph{(1)} \quad $\min \limits_{x\in(x^*-\delta,x^*+\delta)}\frac{\mathrm{d}}{\mathrm{d}x}G'(f(x))=\frac{\mathrm{d}}{\mathrm{d}x}G'(f(x))\mid_x={x^{*-}}$ and $x^*$ is the unique point to attain its minimum.}\\\\
\emph{\emph{(2)} \quad $\min \limits_{x\in(x^*-\delta,x^*+\delta)}\frac{\mathrm{d}}{\mathrm{d}x}G'(f(x))=\frac{\mathrm{d}}{\mathrm{d}x}G'(f(x))\mid_x={x^{*+}}$ and $x^*$ is the unique point to attain its minimum.}\\\\\\
$\mathbf{Case\ 3}$\\\\
Since $G$ is smooth, it can be seen that if $f$ is smooth at the interval $\left[x_l,x_r\right]$, then there exist $x_s\in\left[x_l,x_r\right]$ such that\\\\
$$t=-\frac{1}{\frac{\mathrm{d}}{\mathrm{d}x}G'(f(x))\mid_{x=x_s}}=-\frac{1}{G''(f(x_s))f'(x_s)} \eqno(2.1.3)$$\\\\
One needs to know locally when will the first shock form, that is, the smallest time when locally characteristic intersects. Since $t\ge0$, there is no shock forms when $\frac{\mathrm{d}}{\mathrm{d}x}G'(f(x))\ge0$, that is, when $G'(f(x))$ increases; if $\frac{\mathrm{d}}{\mathrm{d}x}G'(f(x))<0$ then the value of the expression (2.2.2) is positive. We shall define the minimum of (2.2.2) to expect the time when potentially the first shock forms.\\\\
\textbf{Definition\ 3:} \emph{A point $x^*$ is called the \textbf{third kind shock point} if the following conditions are satisfied}\\\\
1 \quad \emph{$f'(x)$ is continuous at $x^*$}\\\\
2 $\quad \frac{\mathrm{d}}{\mathrm{d}x}G'(f(x))\mid_{x=x^*}<0$\\\\
3 \quad \emph{There exist $\delta_0>0$ such that for any $\delta\in(0,\delta_0)$, $\min \limits_{x\in(x^*-\delta,x^*+\delta)}\frac{\mathrm{d}}{\mathrm{d}x}G'(f(x))=\frac{\mathrm{d}}{\mathrm{d}x}G'(f(x))\mid_x={x^*}$} \emph{and} $x^*$ \emph{is the unique point to attain its minimum.}\\\\\\
\textbf{Case 4}\\\\ 
There is another kind of shock point which immediately starts to cross the rarefaction when the shock forms.\\\\ 
\textbf{Definition\ 4:} \emph{A point is call the \textbf{fourth  kind shock point} if the following conditions are satisfied}\\\\
1 \quad \emph{$x^*$ is a jump point of $f(x)$ and $G'(f(x^{*-}))<G'(f(x^{*+}))$}\\\\
2 \quad \emph{$\frac{\mathrm{d}}{\mathrm{d}x}G'(f(x^{*-}))<0$ or $\frac{\mathrm{d}}{\mathrm{d}x}G'(f(x^{*+}))<0$}\\\\
3 \quad \emph{There exist $\delta_0>0$ such that for any $\delta\in(0,\delta_0)$, at least one of the following \emph{(1) (2)} is satisfied}\\\\
\emph{\emph{(1)} \quad $\min \limits_{x\in\left(x^*-\delta,x^*\right]}\frac{\mathrm{d}}{\mathrm{d}x}G'(f(x))=\frac{\mathrm{d}}{\mathrm{d}x}G'(f(x))\mid_x={x^*}$ and $x^*$ is the unique point to attain its minimum.}\\\\
\emph{\emph{(2)} \quad $\min \limits_{x\in\left[x^*,x^*+\delta\right)}\frac{\mathrm{d}}{\mathrm{d}x}G'(f(x))=\frac{\mathrm{d}}{\mathrm{d}x}G'(f(x))\mid_x={x^*}$ and $x^*$ is the unique point to attain its minimum.}\\\\\\
To see the examples and figures of these 4 different kinds of shock points, please go to the section 4.3.\\\\\\\\
\noindent Now from the their definition we have the following propositions revealing the relationship between the break times and their corresponding shock points.\\\\
\textbf{Proposition 1}: \emph{If point $x^*$ is a shock point which is not the fourth kind and denote the break time $t_b$, then}\\\\
\emph{$t_b$=\ \emph{0} \hfill{if $x^*$ is first kind break point.}}\\\\
\emph{$t_b$=\ $\min\{-\frac{1}{\frac{\mathrm{d}}{\mathrm{d}x}G'(f(x))\mid_{x=x^{*-}}},-\frac{1}{\frac{\mathrm{d}}{\mathrm{d}x}G'(f(x))\mid_{x=x^{*+}}}\}$ \hfill{if $x^*$ is second kind break point.}}\\\\
\emph{$t_b$=\ $-\frac{1}{\frac{\mathrm{d}}{\mathrm{d}x}G'(f(x))\mid_{x=x^*}}$ \hfill{if $x^*$ is third kind break point.}}\\\\\\
As for the fourth kind shock point, it has at least one break time, or sometimes two, depending on the condition of its side derivatives.\\\\
\textbf{Proposition 2}: \emph{If point $x^*$ is the fourth kind shock point, then}\\\\
1 \quad $t_{b_1}=-\frac{1}{\frac{\mathrm{d}}{\mathrm{d}x}G'(f(x))\mid_{x=x^{*-}}}$ \emph{if \emph{(1)} of the definition $4$ satisfies and $\frac{\mathrm{d}}{\mathrm{d}x}G'(f(x^{*-}))<0$}\\\\
1 \quad $t_{b_2}=-\frac{1}{\frac{\mathrm{d}}{\mathrm{d}x}G'(f(x))\mid_{x=x^{*+}}}$ \emph{if \emph{(2)} of the definition $4$ satisfies and $\frac{\mathrm{d}}{\mathrm{d}x}G'(f(x^{*+}))<0$}\\\\
Here the proof is omitted since it is basically a very geometrical property. One can imagine them as the level movement of the function, and the property of inflection point.\newpage

\subsection{Shock Equation and Exact Solution of Burger's Equation}
We are now going to establish the equations for any kind of shock that might occur.\\\\\\
\textbf{Case 1}\\\\
From (1.2.5) we know the shock speed is a function of $u_r$ and $u_l$ which are the value of $u$ on the left of the shock and on the right of the shock. Assume the shock coordinate $(\xi(t),t)$ is not on the boundary of rarefaction wave, according to the property of characteristic and the previous section, we must have the following at the shock coordinate $(\xi(t),t)$\\\\
$$u_l=f(x_l) \qquad and \qquad u_r=f(x_r) \qquad for \ \ some \ \ x_l \ <\ x^*\ <\  \ x_r \eqno(2.2.1)$$\\\\
Then we can use (2.1.1) to write down the following equations for 2 characteristic at the point $(\xi(t),t)$.\\\\
$$
\begin{cases}
G'(f(x_l))t=\xi(t)-x_l\\
G'(f(x_r))t=\xi(t)-x_r\\
\end{cases}
\eqno(2.2.2)$$\\\\
We treat $x_l$, $x_r$ both as the functions of t. Take the derivative with respect to t at (2.2.2) we obtain\\\\

$$
\begin{cases}
G''(f(x_l))f'(x_l)\frac{\mathrm{d}x_l}{\mathrm{d}t}t+G'(f(x_l))=\xi'(t)-\frac{\mathrm{d}x_l}{\mathrm{d}t}\\
G''(f(x_r))f'(x_r)\frac{\mathrm{d}x_r}{\mathrm{d}t}t+G'(f(x_r))=\xi'(t)-\frac{\mathrm{d}x_r}{\mathrm{d}t}\\
\end{cases}
\eqno(2.2.3)$$\\\\
Then we substitute $\xi'(t)$ with the jump condition (1.2.5) and use the expression of (2.2.1)\\\\

$$
\begin{cases}
G''(f(x_l))f'(x_l)\frac{\mathrm{d}x_l}{\mathrm{d}t}t+G'(f(x_l))=\frac{G(f(x_r)-G(f(x_l))}{f(x_r)-f(x_l)}-\frac{\mathrm{d}x_l}{\mathrm{d}t}\\
G''(f(x_r))f'(x_r)\frac{\mathrm{d}x_r}{\mathrm{d}t}t+G'(f(x_r))=\frac{G(f(x_r)-G(f(x_l))}{f(x_r)-f(x_l)}-\frac{\mathrm{d}x_r}{\mathrm{d}t}\\
\end{cases}
\eqno(2.2.4)$$\\\\
Move the term $G'(f)$ from left to right, and the $\frac{\mathrm{d}}{\mathrm{d}t}$ from right to left, eventually we obtain

$$
\begin{cases}
\frac{\mathrm{d}x_l}{\mathrm{d}t}=\frac{\frac{G(f(x_r)-G(f(x_l))}{f(x_r)-f(x_l)}-G'(f(x_l))}{1+G''(f(x_l))f'(x_l)t}\\
\\
\frac{\mathrm{d}x_r}{\mathrm{d}t}=\frac{\frac{G(f(x_r)-G(f(x_l))}{f(x_r)-f(x_l)}-G'(f(x_r))}{1+G''(f(x_r))f'(x_r)t}\\
\end{cases}
\eqno(2.2.5)$$\\\\
Thus, we now have a system of ODE. We can also write it in the following form\\\\
$$\frac{\mathrm{d}\mathbf{x}}{\mathrm{d}t}=\mathbf{D}(\mathbf{x},t) \eqno(2.2.6)$$\\\\
where $\mathbf{x}=\begin{bmatrix}x_l\\x_r\end{bmatrix}$  and $\mathbf{D}(\mathbf{x},t)=\begin{bmatrix}\frac{\frac{G(f(x_r)-G(f(x_l))}{f(x_r)-f(x_l)}-G'(f(x_l))}{1+G''(f(x_l))f'(x_l)t} \\\\ \frac{\frac{G(f(x_r)-G(f(x_l))}{f(x_r)-f(x_l)}-G'(f(x_r))}{1+G''(f(x_r))f'(x_r)t}\end{bmatrix}.$\\\\\\\\
We call (2.2.5) or (2.2.6) the \textbf{first kind shock wave equation}.\\\\\\
\textbf{Case 2}\\\\
There exists another case in which shock appear adjacent to the rarefaction, that is, the shock wave go across the region of rarefaction wave. Without loss of generality, we analyze the case in which the shock wave cross the rarefaction wave from the left side. Set the jump point of rarefaction is $x^*$, the following equations come from (1.2.5), (1.3.4), (2.1.1), (2.2.5)\\\\
$$
\begin{cases}
\frac{\mathrm{d}x_l}{\mathrm{d}t}=\frac{\frac{G(u_r)-G(u_l)}{u_r-u_l}-G'(f(x_l))}{1+G''(f(x_l))f'(x_l)t}\\
u_r=(G')^{-1}(\frac{x-x^*}{t})\\
G'(f(x_l))t=x-x_l\\
\end{cases}
$$\\\\
Substitute $x$ with expression of $t$ and $x_l$, then substitute $u_r$ with expression of $t$ and $x_l$, we obtain
$$\frac{\mathrm{d}x_l}{\mathrm{d}t}=\frac{\frac{G((G')^{-1}(G'(f(x_l))+\frac{x_l-x^*}{t}))-G(f(x_l))}{(G')^{-1}(G'(f(x_l))+\frac{x_l-x^*}{t})-f(x_l)}-G'(f(x_l))}{1+G''(f(x_l))f'(x_l)t} \eqno(2.2.7)$$\\\\
We call (2.2.7) the \textbf{second kind shock wave equation}.\\\\\\
\textbf{Case 3}\\\\
Finally there is another case which is somewhat trivial: a rarefaction cross another rarefaction. This indeed might happen, and in this case the shock curve will not related to the any points of initial condition any more -- it is just a curve which can be solve independently. Assume there are 2 rarefaction points $x_1^*$ and $x_2^*$, then the equation is\\\\
$$\frac{\mathrm{d}x}{\mathrm{d}t}=\frac{G((G')^{-1}(\frac{x-x_1^*}{t}))-G((G')^{-1}(\frac{x-x_2^*}{t}))}{(G')^{-1}(\frac{x-x_1^*}{t})-(G')^{-1}(\frac{x-x_2^*}{t})} \eqno(2.2.8)$$\\\\
We call (2.2.8) the \textbf{third kind shock wave equation}.\\\\\\
When $G(u)=\frac{1}{2}u^2$ we call the equation \textbf{Burger's Equation}. We are lucky to see that in this case thees shock equations can be solve analytically.\\\\
\textbf{Property 1:} \emph{If $G(u)=\frac{1}{2}u^2$, that is, for Burger's Equation, the $x_l$ and $x_r$ of \emph{(2.2.5)} satisfies}
$$\frac{1}{2}(f(x_l)+f(x_r))=\frac{\int_{x_l}^{x_r}f(x)\mathrm{d}x+C}{x_r-x_l} \eqno(2.2.9)$$\\\\
\emph{for some constan C. The solution for \emph{(2.2.7)} is}\\\\
$$t=\frac{(x_l-x^*)^2}{2\int_{x^*}^{x_l}(x^*-x)f'(x)\mathrm{d}x+C} \eqno(2.2.10)$$\\\\
\emph{for some constant C. And the solution for \emph{(2.2.8)} is}\\\\
$$x=\frac{x_1^*+x_2^*}{2}+Ct \eqno(2.2.11)$$
\emph{for some constant C.}\\\\
\textbf{\emph{Proof}}: Set $G(u)=\frac{1}{2}u^2$, then (2.2.5) becomes:\\\\
$$
\begin{cases}
\frac{\mathrm{d}x_l}{\mathrm{d}t}=\frac{1}{2}\frac{f(x_r)-f(x_l)}{1+f'(x_l)t}\\
\frac{\mathrm{d}x_r}{\mathrm{d}t}=\frac{1}{2}\frac{f(x_l)-f(x_r)}{1+f'(x_r)t}
\end{cases}
$$\\\\
Divide the first equation by the second equation, we have\\\\
$$\frac{\mathrm{d}x_l}{\mathrm{d}x_r}=-\frac{1+f'(x_r)t}{1+f'(x_l)t}$$\\\\
Substitute $t$ with (2.1.2) it is\\\\
$$\frac{\mathrm{d}x_l}{\mathrm{d}x_r}=-\frac{1-f'(x_r)\frac{x_r-x_l}{f(x_r)-f(x_l)}}{1-f'(x_l)\frac{x_r-x_l}{f(x_r)-f(x_l)}}$$\\\\
which is\\\\
$$\left[f(x_r)-f(x_l)-f'(x_l)(x_r-x_l))\right]\mathrm{d}x_l+\left[f(x_r)-f(x_l)-f'(x_r)(x_r-x_l))\right]\mathrm{d}x_r=0 \eqno(2.2.12)$$\\\\
It is easy to verify the left side of (2.2.12) is a exact differentiation. Intergrate (2.2.12) gives:\\\\
$$(x_r-x_l)\left[f(x_l)+f(x_r)\right]-2\int_{x_l}^{x_r}f(x)\mathrm{d}x=C.$$\\\\
which is exactly (2.2.9).\\\\
As for (2.2.7), Burger's equation reduces it as the following:\\\\
$$\frac{\mathrm{d}x_l}{\mathrm{d}t}=\frac{x_l-x^*}{2t\left[1+f'(x_l)t\right]}$$\\\\
We flip it to make it a Bernoulli equation\\\\
$$\frac{\mathrm{d}t}{\mathrm{d}x_l}=\frac{2}{x_l-x^*}t+\frac{2f'(x_l)}{x_l-x^*}t^2$$\\\\
which can be solved exactly as (2.2.10).\\\\
For the case of the third kind shock equation, now it becoms:\\\\
$$\frac{\mathrm{d}x}{\mathrm{d}t}=\frac{1}{2}(\frac{x-x_1^*}{t}+\frac{x-x_2^*}{t})$$\\\\
Which is\\\\
$$\frac{\mathrm{d}x}{\mathrm{d}t}-\frac{x}{t}=-\frac{x_1^*+x_2^*}{t}$$\\\\
This is a first order linear ODE, which can be solved exactly as (2.2.11) $\Box$.\\\\\\
\textbf{Remark}: We are actually discussing the local properties of the shock point. We call in this way however, a shock point does not necessarily become a shock in various cases. And in the case when at point $x^*$ where $f$ is not smooth, all equations and formulas still work by denoting each term contains $x_r$ as the right side function and the term contains $x_l$ the left side function. One can try to set $f_+(x)$ and $f_-(x)$ and verify that indeeds all the content above are in the same form and they are almost uniform.
\newpage

\section{Initial Point Singularity Analysis of Shock Equation}
In last section we have obtained mathematically very beautiful ODE for shock waves, which means the difficulty of applying characteristic can be solved somehow. At the same time we also have the break point exactly. It seems that we can immediately use the break point as somewhat a initial value of the shock ODE and then solve it. However, one might have realized a very serious problem: in some cases at break point the ODE is singular. This section is to solve this difficulty numerically, and finally build a outline of a new numerical method to solve (1.1.1).
\subsection{Local Differential Property}
For the first kind shock point, the shock forms immediately, and (2.2.6) has no singularity at shock point. We then can directly use it to locally predict the shock. Now let's see what about the rest cases.\\\\
In the case of second kind shock point, from \textbf{proposition 1} we know the break time is the smallest one on the both side. Without loss of generality we say the right side break time is smaller. Then from (2.2.6) $\frac{\mathrm{d}x_l}{\mathrm{d}t}$ is 0 but $\frac{\mathrm{d}x_r}{\mathrm{d}t}$ is like $\frac{0}{0}$. The natural initial point for the equation $(-\frac{1}{G''(f(x^*))f'(x^*)},x^*,x^*)$ is illegal. We have to find another initial point but of course, it is impossible to find a point that is exactly on the solution without solving the equation explicitly. Thankfully, the mathematical techniques help us to numerically find one.\\\\
\textbf{Theorem 1}: \emph{Suppose $x^*$ is second kind shock point, then:}
$$\lim_{t \to t_b}\frac{\mathrm{d}x_l}{\mathrm{d}x_r}=0 \qquad if \qquad \frac{\mathrm{d}}{\mathrm{d}x}G'(f(x))\mid_{x={x^{*+}}}<\frac{\mathrm{d}}{\mathrm{d}x}G'(f(x))\mid_{x={x^{*-}}}$$ \\
$$\lim_{t \to t_b}\frac{\mathrm{d}x_r}{\mathrm{d}x_l}=0 \qquad if \qquad \frac{\mathrm{d}}{\mathrm{d}x}G'(f(x))\mid_{x={x^{*-}}}<\frac{\mathrm{d}}{\mathrm{d}x}G'(f(x))\mid_{x={x^{*+}}}$$ \\\\
\textbf{Proof}: Without loss of generality we assume $\frac{\mathrm{d}}{\mathrm{d}x}G'(f(x))\mid_{x={x^{*+}}}<\frac{\mathrm{d}}{\mathrm{d}x}G'(f(x))\mid_{x={x^{*-}}}$. Devide the two equations of (2.2.6) we obtain\\\\
$$\frac{\mathrm{d}x_l}{\mathrm{d}x_r}=\frac{\frac{G(f(x_r)-G(f(x_l))}{f(x_r)-f(x_l)}-G'(f(x_l))}{\frac{G(f(x_r)-G(f(x_l))}{f(x_r)-f(x_l)}-G'(f(x_r))}\cdot\frac{1+G''(f(x_r))f_+'(x_r)t}{1+G''(f(x_l))f_-'(x_l)t}$$\\\\
Denote the first term\\\\
$$A(x_l,x_r)=\frac{G(f(x_r))-G(f(x_l))-G'(f(x_l))(f(x_r)-f(x_l))}{G(f(x_r))-G(f(x_l))-G'(f(x_r))(f(x_r)-f(x_l))}$$\\\\
and second term\\\\
$$B(x_l,x_r)=\frac{1+G''(f(x_r))f'(x_r)t}{1+G''(f(x_l))f'(x_l)t}$$\\\\
Since $G$ is smooth we use its taylor expansion\\\\
$$G(f(x_r))=G(f(x_l))+\sum_{k=1}^\infty \frac{G^{(k)}(f(x_l))}{k!}(f(x_r)-f(x_l))^k$$\\\\
substitute it into $A(x_l,x_r)$ we have\\\\
$$A(x_l,x_r)=\frac{\sum_{k=2}^\infty \frac{G^{(k)}(f(x_l))}{k!}(f(x_r)-f(x_l))^k}{(G'(f(x_l))-G'(f(x_r)))(f(x_r)-f(x_l))+\sum_{k=2}^\infty \frac{G^{(k)}(f(x_l))}{k!}(f(x_r)-f(x_l))^k}$$\\\\
use the taylor expansion for $G'$\\\\
$$G'(f(x_r))=G'(f(x_l))+\sum_{k=2}^\infty \frac{G^{(k)}(f(x_l))}{(k-1)!}(f(x_r)-f(x_l))^{k-1}$$\\\\
then we obtain\\\\
$$A(x_l,x_r)=\frac{\sum_{k=2}^\infty \frac{G^{(k)}(f(x_l))}{k!}(f(x_r)-f(x_l))^{k-2}}{\sum_{k=2}^\infty (\frac{1}{k!}-\frac{1}{(k-1)!})G^{(k)}(f(x_l))(f(x_r)-f(x_l))^{k-2}}$$\\\\
From the definition of the second kind shock point, we must have $G''(f(x^*))\neq0$. Then take the limit\\\\
$$\lim_{t \to t_b}A(x_l,x_r)=\lim_{f(x_r)-f(x_l) \to 0}A(x_l,x_r)=\frac{\frac{1}{2!}}{\frac{1}{2!}-\frac{1}{(2-1)!}}=-1 \eqno (3.1.1)$$\\\\
which is finite. As for $B(x_l,x_r)$, since $\frac{\mathrm{d}}{\mathrm{d}x}G'(f(x))\mid_{x={x^{*+}}}<\frac{\mathrm{d}}{\mathrm{d}x}G'(f(x))\mid_{x={x^{*-}}}$, we have $t_b$=$-\frac{1}{G''(f(x^*))f'(x^{*-})}$. When $t=t_b$ its numerator is zero but denominator is non-zero. Therefore\\\\
$$\lim_{t \to t_b}\frac{\mathrm{d}x_l}{\mathrm{d}x_r}=\lim_{t \to t_b}\left[A(x_l,x_r)\cdot B(x_l,x_r)\right]=-1\cdot0=0$$\\\\
The result for the case when $\frac{\mathrm{d}}{\mathrm{d}x}G'(f(x))\mid_{x={x^{*-}}}<\frac{\mathrm{d}}{\mathrm{d}x}G'(f(x))\mid_{x={x^{*+}}}$ immediately comes from the proof above and the symmetry of $x_l$ and $x_r$ of the equation (2.2.6) $\quad \Box.$\\\\
Now it is time to consider the last case -- what is the local differential property of the thrid kind shock point? Naively from intuition the symmetry gives $\lim_{t \to t_b}\frac{\mathrm{d}x_l}{\mathrm{d}x_r}=-1$. Indeed, this is only guaranteed if $f$ is smooth at $x^*$ as we will see in a short time, and there are many possibilities for this case. We shall prove a much stronger result but before we go into the detailed proof, let us set some lemmas.\\\\
For simplicity, we denote $h(x)\triangleq G'(f(x))$, $h_-^{(k)}=\frac{\mathrm{d}^{k}h}{\mathrm{d}x^k}\mid_{x=x^{*-}}$ and $h_+^{(k)}=\frac{\mathrm{d}^{k}h}{\mathrm{d}x^k}\mid_{x=x^{*+}}$.\\\\
\textbf{Lemma 1}: \emph{Suppose $x^*$ is third kind shock point. If there exists $k_l$ and $k_r$ that are the smallest integers greater equal than $2$ such that $h_-^{(k_l)}$ and $h_+^{(k_r)}$ does not vanish. Then we have}\\\\
$$(-1)^{k_l-1}h_-^{(k_l)}>0 \qquad and \qquad h_+^{(k_r)}>0$$\\\\
\textbf{Proof:} From the definition of the third kind shock point, there exists $\delta_l>0$ and $\delta_r>0$ such that $h_-^{(2)}(x)<0$ at $(x^*-\delta_l,x^*)$ and $h_+^{(2)}(x)>0$ at $(x^*+\delta_r,x^*)$. If $h_-^{(2)}<0$ then $k_l=2$ and $(-1)^{k_l-1}h_-^{(k_l)}>0$ is true. If $h_-^{(2)}=0$ and $h_-^{(2)}\ne0$, to ensure $h_-^{(2)}(x)<0$, there must be a left neighbourhood of $x^*$ such that $h_-^{(2)}$ is increasing, which means $h_-^{(3)}>0$. If $h_-^{(3)}=0$, just following this process and we can see it is indeed $(-1)^{k_l-1}h_-^{(k_l)}>0$. The same technique can be applied to get $h_+^{(k_r)}>0$ $\quad \Box$.\\\\
\textbf{Theorem 2}: \emph{Suppose $x^*$ is third kind shock point, and $k\ge2$ the smallest integer such that at least one of the $h_-^{(k)}$ and $h_+^{(k)}$ does not vanish, then we have the following conclusions}\\\\
1 \quad$\lim_{t \to t_b}\frac{\mathrm{d}x_l}{\mathrm{d}x_r}=0$ if $h_+^{(k)}=0$\\\\
2 \quad$\lim_{t \to t_b}\frac{\mathrm{d}x_l}{\mathrm{d}x_r}=-\infty$ if $h_-^{(k)}=0$\\\\
3 \quad$\lim_{t \to t_b}\frac{\mathrm{d}x_l}{\mathrm{d}x_r}$ is the negative root of the following polynomial $H(s)$ if none of the $h_+^{(k)}$ and $h_-^{(k)}$ vanishes\\\\
$$H(s)=s^{k+1}-(\frac{k+1}{k-1})s^k+(\frac{k+1}{k-1})ps-p$$\\\emph{where $p=\frac{h_+^{(k)}}{h_-^{(k)}}$.}\\\\
\textbf{Proof}: For $A(x_l,x_r)$ in theorem 1 there is no difference between in its case and in this case. From (3.1.1) we know the limit of $A(x_l,x_r)$ is still $-1$. However for $B(x_l,x_r)$ there is no way to evaluate itself independently. But first at least we can do taylor expasion and see what is going for $B(x_l,x_r)$.\\\\
Substitute $t$ with (2.1.2) into $B(x_l,x_r)$ and use the denotation $h=G'(f)$ we have\\\\
$$B(x_l,x_r)=\frac{h(x_r)-h(x_l)-h'(x_r)(x_r-x_l)}{h(x_r)-h(x_l)-h'(x_l)(x_r-x_l)}$$\\\\
For simplicity we denote $h=h(x^*)$, $h'=h'(x^*)$, $\triangle x_l=(x_l-x^*)$ and  $\triangle x_r=(x_r-x^*)$. Note that $x_r-x_l=\triangle x_r-\triangle x_l$. Use the fact that\\\\
$$h(x_l)=h+h'\triangle x_l+\sum_{n=2}^\infty \frac{h_-^{(n)}}{n!}\triangle x_l^n$$ $$h'(x_l)=h'+\sum_{n=2}^\infty \frac{h_-^{(n)}}{(n-1)!}\triangle x_l^{n-1}$$ $$h(x_r)=h+h'\triangle x_r+\sum_{n=2}^\infty \frac{h_+^{(n)}}{n!}\triangle x_r^n$$ $$h'(x_r)=h'+\sum_{n=2}^\infty \frac{h_+^{(n)}}{(n-1)!}\triangle x_r^{n-1}$$\\
we can get\\\\
$$B(x_l,x_r)=\frac{\sum_{n=2}^\infty\frac{1}{n!}(h_{+}^{(n)}\triangle x_r^n-h_{-}^{(n)}\triangle x_l^n)-(\triangle x_r-\triangle x_l)\sum_{n=2}^\infty\frac{1}{(n-1)!}h_{+}^{(n)}\triangle x_r^{n-1}}{\sum_{n=2}^\infty\frac{1}{n!}(h_{+}^{(n)}\triangle x_r^n-h_{-}^{(n)}\triangle x_l^n)-(\triangle x_r-\triangle x_l)\sum_{n=2}^\infty\frac{1}{(n-1)!}h_{-}^{(n)}\triangle x_l^{n-1}}$$\\\\
Let $k$ be the smaller integer such that at least one of the $h_{-}^{(k)}$ and $h_{+}^{(k)}$ does not vanish. Then\\\\
$$B(x_l,x_r)=\frac{\sum_{n=k}^\infty\frac{1}{n!}(h_{+}^{(n)}\triangle x_r^n-h_{-}^{(n)}\triangle x_l^n)-(\triangle x_r-\triangle x_l)\sum_{n=k}^\infty\frac{1}{(n-1)!}h_{+}^{(n)}\triangle x_r^{n-1}}{\sum_{n=k}^\infty\frac{1}{n!}(h_{+}^{(n)}\triangle x_r^n-h_{-}^{(n)}\triangle x_l^n)-(\triangle x_r-\triangle x_l)\sum_{n=k}^\infty\frac{1}{(n-1)!}h_{-}^{(n)}\triangle x_l^{n-1}}$$\\\\
Devide both numerator and denominator by $\triangle x_r^k$\\\\
$$B(x_l,x_r)=\frac{\frac{1}{k!}(h_+^{(k)}-h_-^{(k)}(\frac{\triangle x_l}{\triangle x_r})^k)-\frac{1}{(k-1)!}h_+^{(k)}(1-\frac{\triangle x_l}{\triangle x_r})+O(\triangle x_r)}{\frac{1}{k!}(h_+^{(k)}-h_-^{(k)}(\frac{\triangle x_l}{\triangle x_r})^k)-\frac{1}{(k-1)!}h_-^{(k)}(1-\frac{\triangle x_l}{\triangle x_r})(\frac{\triangle x_l}{\triangle x_r})^{k-1}+O(\triangle x_l)}$$\\\\
Now we can write\\\\
$$\frac{\mathrm{d}x_r}{\mathrm{d}x_l}=A(x_l,x_r)\cdot B(x_l,x_r)=A(x_l,x_r)\cdot\frac{(h_+^{(k)}-h_-^{(k)}(\frac{\triangle x_l}{\triangle x_r})^k)-kh_+^{(k)}(1-\frac{\triangle x_l}{\triangle x_r})+O(\triangle x_r)}{(h_+^{(k)}-h_-^{(k)}(\frac{\triangle x_l}{\triangle x_r})^k)-kh_-^{(k)}(1-\frac{\triangle x_l}{\triangle x_r})(\frac{\triangle x_l}{\triangle x_r})^{k-1}+O(\triangle x_l)}$$\\\\
Now take the limit with $t\to t_b$. Note that $\lim_{t\to t_b}\triangle x_l=\lim_{t\to t_b}\triangle x_r=0$ and $\lim_{t\to t_b}\frac{\triangle x_l}{\triangle x_r}=\lim_{t\to t_b}\frac{\mathrm{d}x_r}{\mathrm{d}x_l}$. Set $\lim_{t\to t_b}\frac{\mathrm{d}x_r}{\mathrm{d}x_l}=s$, use (3.1.1) we obtain\\\\
$$s=-\frac{(h_+^{(k)}-h_-^{(k)}s^k)-kh_+^{(k)}(1-s)}{(h_+^{(k)}-h_-^{(k)}s^k)-kh_-^{(k)}(1-s)s^{k-1}} \eqno (3.1.2)$$\\\\
We first suppose none of the $h_-^{(k)}$ and $h_+^{(k)}$ vanishes. Then we can let $p=\frac{h_+^{(k)}}{h_-^{(k)}}$, and one can see that (3.1.2) is exactly the polynomial equation:\\\\
$$s^{k+1}-(\frac{k+1}{k-1})s^k+(\frac{k+1}{k-1})ps-p=0 \eqno (3.1.3)$$\\\\
In the case of $h_+^{(k)}$ vanishes, (3.1.2) still holds and there is only a non-positive root which is 0, and this implies $\lim_{t \to t_b}\frac{\mathrm{d}x_l}{\mathrm{d}x_r}=0$; if $h_-^{(k)}=0$, (3.1.2) only has one root s=$\frac{k-1}{k+1}$ which is positive. Of course according to the symmetry of the equation one can in another way, prove $\lim_{t \to t_b}\frac{\mathrm{d}x_l}{\mathrm{d}x_r}=-\infty$ by the flipping it as $\lim_{t \to t_b}\frac{\mathrm{d}x_r}{\mathrm{d}x_l}=0$. But in fact we can just treat it as the case as $h_-^{(k)}\to0$ and see what happens to the (3.1.3).\\\\
\textbf{Lemma 2}: \emph{The polynomial function $H(s)$ has a unique negative root if $p\ne0$. Moreover, if the root is $r(p)$, then $(-1)^k\cdot\frac{\mathrm{d}r}{\mathrm{d}p}>0$.}\\\\
\textbf{Proof:} Take the first and second derivative of $H(s)$ we have\\\\
$$H(s)=s^{k+1}-(\frac{k+1}{k-1})s^k+(\frac{k+1}{k-1})ps-p$$\\
$$H'(s)=(k+1)(s^k-\frac{k}{k-1}s^{k-1}+\frac{p}{k-1})$$\\
$$H''(s)=k(k+1)s^{k-2}(s-1) \eqno(3.1.4)$$\\\\
We first prove the case when $k$ is even. (3.1.4) shows in this case $H'(s)$ decreases when $s<0$. From lemma 1 we know now $p<0$, so when $s<0$ $H'(s)$ decreases from $+\infty$ and cross the negative axis and arrive $(0,\frac{k+1}{k-1}p)$ which is on the lower-halfplane. Along with $H(0)=-p>0$, this means $H(s)$ must have a unique negative root; if $k$ is odd, now $p>0$, (3.1.4) tells $H'(s)$ increases when $s<0$ and arrive $(0,\frac{k+1}{k-1}p)$ which is on the upper-halfplane. With now $H(0)=-p<0$ we know $H(s)$ still have a unique negative root.\\\\
The root $r(p)$ satisfies $H(r)=0$. Differential this equation with respect to $r$ and $p$ we have\\\\
$$H'(r)dr+(\frac{k+1}{k-1}r-1)dp=0$$\\\\
Which is also\\\\
$$\frac{\mathrm{d}r}{\mathrm{d}p}=-\frac{\frac{k+1}{k-1}r-1}{H'(r)}$$\\\\
From the proof of the uniqueness of the negative root, $H'(r)>0$ when $k$ is even and $H'(r)<0$ when $k$ is odd. The negative root $r(p)$ always makes the numerator $\frac{k+1}{k-1}r-1<0$. This implies $(-1)^k\cdot\frac{\mathrm{d}r}{\mathrm{d}p}>0$.\\\\
Now we are able to complete the proof of theorem 2 more elegantly. In either the case $k$ is even or $k$ is odd, either $h_-^{(k)}=0$ or $h_+^{(k)}=0$, we treat them as the limit taken from the correct side and find that the negative root goes in the way that is uniform to the conclusions of theorem 2 \quad $\Box$.\\\\
\textbf{Corollary 1:} \emph{If f is smooth at $x^*$, then $\lim_{t \to t_b}\frac{\mathrm{d}x_l}{\mathrm{d}x_r}=-1$.}\\\\
\textbf{Proof:} We can actually say it in a slightly stronger way. One can substitute $s=-1$ and see that it is indeed the root of $H(s)$ in this case of $p=\frac{h_+^{(k)}}{h_-^{(k)}}=(-1)^{k-1}$. According to the lemma 2 this negative root is unique, so $\lim_{t \to t_b}\frac{\mathrm{d}x_l}{\mathrm{d}x_r}$ cannot be the value other than $-1$ \quad $\Box$.\\\\
\newpage
\subsection{Local Asymptotic Differentiation}
In the last subsection, we have obtained some properties about the shock equation near the local shock point, however it does not mean we are already able to choose a legal initial point of the equation. The reason is simple: perhaps we still need to think more about the case when $\frac{\mathrm{d}x_l}{\mathrm{d}x_r}$ vanishes or blows up. In the numerical simulation, for example, an unstable shock curve might occur if one simply just use $(\frac{\delta}{h(x^*)-h(x^*+\delta)},x^*,x^*+\delta)$ as initial data in the case of $\frac{\mathrm{d}x_l}{\mathrm{d}x_r}=0$. Asymptotic results stronger than $\frac{\mathrm{d}x_l}{\mathrm{d}x_r}=0$ or $\frac{\mathrm{d}x_r}{\mathrm{d}x_l}=0$ are necessary.\\\\
\textbf{Theorem 3:} \emph{Suppose $x^*$ is second kind shock point. Let $k_r$, $k_l$ be that of the lemma 1, then}\\\\
1 \quad $\lim_{t \to t_b}\frac{\triangle x_l}{(\triangle x_r)^{k_r}}=\frac{(k_r-1)h_+^{(k_r)}}{(k_r+1)!(h_+'-h_-')} \quad $ if \quad $h_+'<h_-'$\\\\
2 \quad $\lim_{t \to t_b}\frac{\triangle x_l}{(\triangle x_r)^{k_r}}=\frac{(k_l-1)h_-^{(k_l)}}{(k_l+1)!(h_-'-h_+')}\quad $ if \quad $h_-'<h_+'$\\\\
\textbf{Proof:} We only need to prove the first conclusion as the second one comes immediately from the symmetry of the equations. The following can be obtained if one use the same techniques of that used in the proof of theorem 2\\\\
$$B(x_l,x_r)=\frac{(h_+'-h_-')\triangle x_l+\sum_{n=k_r}^\infty\frac{h_+^{(n)}}{n!}\triangle x_r^n-\sum_{n=k_l}^\infty\frac{h_-^{(n)}}{n!}\triangle x_l^n-(\triangle x_r-\triangle x_l)\sum_{n=k_r}^\infty\frac{h_+^{(n)}}{(n-1)!}\triangle x_r^{n-1}}{(h_+'-h_-')\triangle x_r+\sum_{n=k_r}^\infty\frac{h_+^{(n)}}{n!}\triangle x_r^n-\sum_{n=k_l}^\infty\frac{h_-^{(n)}}{n!}\triangle x_l^n-(\triangle x_r-\triangle x_l)\sum_{n=k_l}^\infty\frac{h_-^{(n)}}{(n-1)!}\triangle x_l^{n-1}}$$\\\\
$$=\frac{(h_+'-h_-')\triangle x_l+\frac{h_+^{(k_r)}}{k_r!}\triangle x_r^{k_r}-\frac{h_-^{(k_l)}}{k_l!}\triangle x_l^{k_l}-(\triangle x_r-\triangle x_l)\frac{h_+^{(k_r)}}{(k_r-1)!}\triangle x_r^{k_r-1}+O(\triangle x_r^{k_r+1})++O(\triangle x_l^{k_l+1})}{(h_+'-h_-')\triangle x_r+\frac{h_+^{(k_r)}}{k_r!}\triangle x_r^{k_r}-\frac{h_-^{(k_l)}}{k_l!}\triangle x_l^{k_l}-(\triangle x_r-\triangle x_l)\frac{h_-^{(k_l)}}{(k_l-1)!}\triangle x_l^{k_l-1}+O(\triangle x_r^{k_r+1})+O(\triangle x_l^{k_l+1})}$$\\\\
We ignore the infinitesimal terms and use asymptotic sign. Along with (3.1.1) we have\\\\
$$\frac{\mathrm{d}x_l}{\mathrm{d}x_r}\sim-\frac{(h_+'-h_-')\triangle x_l+\frac{h_+^{(k_r)}}{k_r!}\triangle x_r^{k_r}-\frac{h_-^{(k_l)}}{k_l!}\triangle x_l^{k_l}-(\triangle x_r-\triangle x_l)\frac{h_+^{(k_r)}}{(k_r-1)!}\triangle x_r^{k_r-1}}{(h_+'-h_-')\triangle x_r+\frac{h_+^{(k_r)}}{k_r!}\triangle x_r^{k_r}-\frac{h_-^{(k_l)}}{k_l!}\triangle x_l^{k_l}-(\triangle x_r-\triangle x_l)\frac{h_-^{(k_l)}}{(k_l-1)!}\triangle x_l^{k_l-1}} \eqno(3.2.1)$$\\\\
Note that $\mathrm{d}x_l=\mathrm{d}(x_l-x^*)=\mathrm{d}\triangle x_l$ and $\mathrm{d}x_l=\mathrm{d}(x_l-x^*)=\mathrm{d}\triangle x_l$. We know that in this case $\frac{\mathrm{d}x_l}{\mathrm{d}x_r}\sim0$. Assume $\triangle x_l\sim p\triangle x_r^q$ from some constant $p$ and $q>1$. Then we have\\\\
$$\frac{\triangle x_l}{\triangle x_r}\sim\frac{\mathrm{d}x_l}{\mathrm{d}x_r}=\frac{\mathrm{d}\triangle x_l}{\mathrm{d}\triangle x_r}\sim\frac{\mathrm{d}(p\triangle x_r^q)}{\mathrm{d}\triangle x_r}\sim pq\triangle x_r^{q-1}$$\\\\
Then (3.2.1) can be written as:\\\\
$$pq(\triangle x_r)^{q-1}\sim-\frac{(h_+'-h_-')p\triangle x_r^q+\frac{h_+^{(k_r)}}{k_r!}\triangle x_r^{k_r}-\frac{h_-^{(k_l)}}{k_l!}p^{k_l}\triangle x_r^{qk_l}-(\triangle x_r-p\triangle x_r^q)\frac{h_+^{(k_r)}}{(k_r-1)!}\triangle x_r^{k_r-1}}{(h_+'-h_-')\triangle x_r+\frac{h_+^{(k_r)}}{k_r!}\triangle x_r^{k_r}-\frac{h_-^{(k_l)}}{k_l!}p^{k_l}\triangle x_r^{qk_l}-(\triangle x_r-p\triangle x_r^q)\frac{h_-^{(k_l)}}{(k_l-1)!}p^{k_l-1}\triangle x_r^{q(k_l-1)}}$$\\\\
Since now $k_r>2$, $k_l>2$, $q>1$, we drop those terms which must not be the lowest order and obtain
$$pq(\triangle x_r)^{q-1}\sim\frac{-(h_+'-h_-')p\triangle x_r^q+\left[\frac{1}{(k_r-1)!}-\frac{1}{k_r!}\right]h_+^{(k_r)}\triangle x_r^{k_r}}{(h_+'-h_-')\triangle x_r}$$\\\\
which is\\\\
$$pq\sim\frac{-(h_+'-h_-')p+\left[\frac{1}{(k_r-1)!}-\frac{1}{k_r!}\right]h_+^{(k_r)}\triangle x_r^{k_r-q}}{h_+'-h_-'}$$\\\\
Since $q>1$, the only possible case is $q=k_r$\\\\
$$pk_r\sim\frac{-(h_+'-h_-')p+\left[\frac{1}{(k_r-1)!}-\frac{1}{k_r!}\right]h_+^{(k_r)}}{h_+'-h_-'}$$\\\\
Solve the $p$ and we obtain\\\\
$$p=\frac{(k_r-1)h_+^{(k_r)}}{(k_r+1)!(h_+'-h_-')}$$\\\\
Thus we have\\\\
$$\triangle x_l\sim\frac{(k_r-1)h_+^{(k_r)}}{(k_r+1)!(h_+'-h_-')}\triangle x_r^{k_r}$$\\\\
For the case when $h_-'<h_+'$, from the symmetry of get the similar result:\\\\
$$\triangle x_r\sim\frac{(k_l-1)h_+^{(k_l)}}{(k_l+1)!(h_-'-h_+')}\triangle x_l^{k_l} \quad \Box.$$\\\\
\textbf{Theorem 4:} \emph{Suppose $x^*$ is third kind shock point, then}\\\\
1 \quad $\triangle x_l\sim-\left|\left[\frac{(1-k_r)k_l!h_+^{(k_r)}}{(k_r+1)!h_-^{(k_l)}}\right]^{\frac{1}{k_l}}\right|\triangle x_r^{\frac{k_r}{k_l}}$ \quad if \quad $k_l<k_r$\\\\
2 \quad $\triangle x_r\sim\left[(-1)^{k_l}\frac{(1-k_l)k_r!h_-^{(k_l)}}{(k_l+1)!h_+^{(k_r)}}\right]^{\frac{1}{k_r}}(-\triangle x_l)^{\frac{k_l}{k_r}}$ \quad if \quad $k_r<k_l$\\\\\\\\
\textbf{Proof:} We first prove the first case. In this case we have
$$B(x_l,x_r)=\frac{\frac{h_+^{(k_r)}}{k_r!}\triangle x_r^{k_r}-\frac{h_-^{(k_l)}}{k_l!}\triangle x_l^{k_l}-(\triangle x_r-\triangle x_l)\frac{h_+^{(k_r)}}{(k_r-1)!}\triangle x_r^{k_r-1}+O(\triangle x_r^{k_r+1})++O(\triangle x_l^{k_l+1})}{\frac{h_+^{(k_r)}}{k_r!}\triangle x_r^{k_r}-\frac{h_-^{(k_l)}}{k_l!}\triangle x_l^{k_l}-(\triangle x_r-\triangle x_l)\frac{h_-^{(k_l)}}{(k_l-1)!}\triangle x_l^{k_l-1}+O(\triangle x_r^{k_r+1})+O(\triangle x_l^{k_l+1})}$$\\\\
Again, drop the higher order terms and suppose $\triangle x_l\sim p\triangle x_r^q$, then we have the asymptotic relation\\\\
$$pq(\triangle x_r)^{q-1}\sim-\frac{\frac{h_+^{(k_r)}}{k_r!}\triangle x_r^{k_r}-\frac{h_-^{(k_l)}}{k_l!}p^{k_l}\triangle x_r^{qk_l}-(\triangle x_r-p\triangle x_r^q)\frac{h_+^{(k_r)}}{(k_r-1)!}\triangle x_r^{k_r-1}}{\frac{h_+^{(k_r)}}{k_r!}\triangle x_r^{k_r}-\frac{h_-^{(k_l)}}{k_l!}p^{k_l}\triangle x_r^{qk_l}-(\triangle x_r-p\triangle x_r^q)\frac{h_-^{(k_l)}}{(k_l-1)!}p^{k_l-1}\triangle x_r^{q(k_l-1)}}$$\\\\
Since $q>1$, drop the terms which must not be the lowest order\\\\
$$pq(\triangle x_r)^{q-1}\sim-\frac{\left[\frac{1}{k_r!}-\frac{1}{(k_r-1)!}\right]h_+^{(k_r)}\triangle x_r^{k_r}-\frac{h_-^{(k_l)}}{k_l!}p^{k_l}\triangle x_r^{qk_l}}{\frac{h_+^{(k_r)}}{k_r!}\triangle x_r^{k_r}-\frac{h_-^{(k_l)}}{(k_l-1)!}p^{k_l-1}\triangle x_r^{qk_l-q+1}}$$\\\\
Note that with $k_l\ge2$ and $k_r\ge2$, one can verify that this is possible only when $k_r=qk_l$. Then we have\\\\
$$p\frac{k_r}{k_l}(\triangle x_r)^{\frac{k_r}{k_l}-1}\sim\frac{\left[\frac{1}{k_r!}-\frac{1}{(k_r-1)!}\right]h_+^{(k_r)}\triangle x_r^{k_r}-\frac{h_-^{(k_l)}}{k_l!}p^{k_l}\triangle x_r^{k_r}}{\frac{h_-^{(k_l)}}{(k_l-1)!}p^{k_l-1}\triangle x_r^{k_r-\frac{k_r}{k_l}+1}}$$\\\\
We are now able to cancle all $\triangle x_r$. Solve with $p$ and note that $p$ must be negative, we obtain\\\\
$$p=-\left|\left[\frac{(1-k_r)k_l!h_+^{(k_r)}}{(k_r+1)!h_-^{(k_l)}}\right]^{\frac{1}{k_l}}\right|$$\\\\
This complete the proof of the first case. For the second case, to avoid the root of negative value $\triangle x_l$, we assume $\triangle x_r\sim p(-\triangle x_l)^q$ for $p>0$ and $q>1$. Then the asymptotic relation becomes\\\\
$$-pq(-\triangle x_l)^{q-1}\sim-\frac{\frac{h_-^{(k_l)}}{k_l!}\triangle x_l^{k_l}-\frac{h_+^{(k_r)}}{k_r!}p^{k_r}(-\triangle x_l)^{qk_r}-(\triangle x_l-p(-\triangle x_l)^q)\frac{h_-^{(k_l)}}{(k_l-1)!}\triangle x_l^{k_l-1}}{\frac{h_-^{(k_l)}}{k_l!}\triangle x_l^{k_l}-\frac{h_+^{(k_r)}}{k_r!}p^{k_r}(-\triangle x_l)^{qk_r}-(\triangle x_l-p(-\triangle x_l)^q)\frac{h_+^{(k_r)}}{(k_r-1)!}p^{k_r-1}(-\triangle x_l)^{q(k_r-1)}}$$\\\\
There is no different for the result of $k_l=qk_r$. To be convenient we write it in this way\\\\
$$-p\frac{k_l}{k_r}(-\triangle x_l)^{\frac{k_l}{k_r}-1}\sim-\frac{(-1)^{k_l}\left[\frac{1}{k_l!}-\frac{1}{(k_l-1)!}\right]h_-^{(k_l)}(-\triangle x_l)^{k_l}-\frac{h_+^{(k_r)}}{k_r!}p^{k_r}(-\triangle x_l)^{k_l}}{\frac{h_+^{(k_r)}}{(k_r-1)!}p^{k_r-1}(-\triangle x_l)^{k_l-\frac{k_l}{k_r}+1}}$$\\\\
So now we are able to cancle all $(-\triangle x_l)$, and obtain\\\\
$$p\frac{k_l}{k_r}\sim\frac{(-1)^{k_l}\left[\frac{1}{k_l!}-\frac{1}{(k_l-1)!}\right]h_-^{(k_l)}-\frac{h_+^{(k_r)}}{k_r!}p^{k_r}}{\frac{h_+^{(k_r)}}{(k_r-1)!}p^{k_r-1}}$$\\\\
Which is 
$$p=\left[(-1)^{k_l}\frac{(1-k_l)k_r!h_-^{(k_l)}}{(k_l+1)!h_+^{(k_r)}}\right]^{\frac{1}{k_r}}$$\\\\
From the lemma 1 one can quickly verify that in either of the cases the sign of $p$ is uniform to the sign of $\triangle x_l$ and $\triangle x_r$ \quad $\Box$.\\\\
\textbf{Theorem 5} \emph{Suppose $x^*$ is fourth kind shock point, then at a small neighbour of $x^*$ the following condition satisfies}\\\\
$$t\sim\frac{1}{-h'(x^*)-\frac{2h^{(k)}(x^*)}{(k-1)!(k+1)}\triangle x^{k-1}}$$\\\\
\textbf{Proof:} Since in such a case as we shall see at the following there is no need to distinguish $x_l$ and $x_r$. We denote both as $x$ (Here $x$ is somewaht still $x_l$ or $x_r$ but not the coordinate of the shock curve $\xi(t)$), and discard any sign that implies the side we are discussing. Recall the second kind shock equation (2.2.7), for simplicity we first denote $u=(G')^{-1}(G'(f(x))+\frac{x-x^*}{t})$ as the solution at the boundary side of rarefaction. Then (2.2.7) can be written as\\\\
$$\frac{\mathrm{d}x}{\mathrm{d}t}=\frac{\frac{G(u)-G(f(x))}{u-f(x)}-G'(f(x))}{1+h'(x)t}$$\\\\
Use the following taylor expansion of $G$\\\\
$$G(u)=G(f(x))+\sum_{m=1}^{\infty}\frac{G^{(m)}(f(x))}{m!}(u-f(x))^m$$\\\\
$$G'(u)=G'(f(x))+\sum_{m=1}^{\infty}\frac{G^{(m+1)}(f(x))}{m!}(u-f(x))^m$$\\\\
We have\\\\
$$\frac{\mathrm{d}x}{\mathrm{d}t}=\frac{\sum_{m=2}^{\infty}\frac{G^{(m)}(f(x))}{m!}(u-f(x))^{m-1}}{1+h'(x)t} \eqno(3.2.2)$$\\\\
Now use the taylor expansion of $(G')^{-1}$\\\\
$$u=(G')^{-1}(G'(f(x))+\frac{x-x^*}{t})=f(x)+\sum_{n=1}^{\infty}\frac{(G'^{(-1)})^{(n)}\mid_{G'(f(x))}}{n!}(\frac{x-x^*}{t})^n$$\\\\
(3.2.2) now becomes\\\\
$$\frac{\mathrm{d}x}{\mathrm{d}t}=\frac{\sum_{m=2}^{\infty}\frac{G^{(m)}(f(x))}{m!}\left[\sum_{n=1}^{\infty}\frac{(G'^{(-1)})^{(n)}\mid_{G'(f(x))}}{n!}(\frac{x-x^*}{t})^n\right]^{m-1}}{1+h'(x)t}$$\\\\
One can verify that the coefficient of $\frac{x-x^*}{t}$ is $\frac{1}{2}$, and the coefficient for the higher order terms of $\frac{x-x^*}{t}$ is related to $m$ and $n$. We write it in this form\\\\
$$\frac{\mathrm{d}x}{\mathrm{d}t}=\frac{x-x^*+\frac{2(x-x^*)^2}{t}P(x,\frac{x-x^*}{t})}{2t\left[1+h'(x)t\right]} \eqno(3.2.3)$$\\\\
where $P(x,\frac{x-x^*}{t})$ is the terms of all higher order of $\frac{x-x^*}{t}$ and is bounded when $x\to x^*$. By flipping (3.2.3) one can verify that the following\\\\
$$\frac{\mathrm{d}t}{\mathrm{d}x}=\frac{2t\left[1+h'(x)t\right]   }{x-x^*+\frac{2(x-x^*)^2}{t}P(x,\frac{x-x^*}{t})}$$\\\\
$$=\frac{2}{x-x^*}t+\frac{2h'(x)}{x-x^*}t^2-\left[1+h'(x)t\right]\frac{4tP(x,\frac{x-x^*}{t})}{t+2(x-x^*)P(x,\frac{x-x^*}{t})} \eqno(3.2.4)$$\\\\
We can solve it somehow as Bernoulli equation. Divide $t^2$ both sides and let $z=\frac{1}{t}$, then (3.2.4) is\\\\
$$\frac{\mathrm{d}z}{\mathrm{d}x}+\frac{2}{x-x^*}z=-2\frac{h'(x)}{x-x^*}+z\left[1+h'(x)t\right]\frac{4tP(x,\frac{x-x^*}{t})}{t+2(x-x^*)P(x,\frac{x-x^*}{t})} \eqno(3.2.5)$$
Since $t_b=-\frac{1}{h'(x^*)}$ which is not zero, $\frac{4tP(x,\frac{x-x^*}{t})}{t+2(x-x^*)P(x,\frac{x-x^*}{t})}$ is bounded or even vanishes, and $1+h'(x)t$ vanishes near $t_b$. Assume when $t\to t_b$ satisfies\\$z\left[1+h'(x)t\right]\sim p_q(x-x^*)^q$ and $\frac{4tP(x,\frac{x-x^*}{t})}{t+2(x-x^*)P(x,\frac{x-x^*}{t})}\sim p_r(x-x^*)^r$ for some constant $p_q$, $p_r$, and $q>0$, $r\ge0$. Then (3.2.5) is\\\\
$$\frac{\mathrm{d}z}{\mathrm{d}x}+\frac{2}{x-x^*}z\sim-2\frac{h'(x)}{x-x^*}+p_qp_r(x-x^*)^{q+r} \eqno(3.2.6)$$\\\\
Multiply both sides of (3.2.6) with $(x-x^*)^2$ then we have\\\\
$$\left[(x-x^*)^2z\right]'\sim-2h'(x)(x-x^*)+p_qp_r(x-x^*)^{q+r+2}$$\\\\
Which is\\\\
$$z\sim\frac{-2\int_{x^*}^{x}h'(s)(s-x^*)ds}{(x-x^*)^2}+\frac{p_qp_r}{q+r+3}(x-x^*)^{q+r+1} \eqno(3.2.7)$$\\\\
Use the taylor expansion of $h'(s)$\\\\
$$h'(x)=\sum_{n=1}^{\infty}\frac{h^{(n)}(x^*)}{(n-1)!}(x-x^*)^{n-1}$$\\\\
Then (3.2.7) is\\\\
$$z\sim-h'(x^*)-\sum_{n=2}^{\infty}\frac{2h^{(n)}(x^*)}{(n-1)!(n+1)}(x-x^*)^{n-1}+\frac{p_qp_r}{q+r+3}(x-x^*)^{q+r+1}$$\\\\
Note that $z+h'(x^*)=z\left[1+h'(x^*)t\right]\sim p_q(x-x^*)^q$, we must have\\\\
$$\frac{-2h^{(k)}(x^*)}{(k-1)!(k+1)}(x-x^*)^{k-1}+\frac{p_qp_r}{q+r+3}(x-x^*)^{q+r+1}\sim p_q(x-x^*)^q$$\\\\
Note that $r\ge0$, we obtain\\\\
$$\frac{-2h^{(k)}(x^*)}{(k-1)!(k+1)}(x-x^*)^{k-1}\sim p_q(x-x^*)^q$$\\\\
Thus $q=k-1$ and $p_q=\frac{-2h^{(k)}(x^*)}{(k-1)!(k+1)}$. Therefore\\\\
$$z\left[1+h'(x)t\right]=\frac{1}{t}+h'(x)\sim\frac{-2h^{(k)}(x^*)}{(k-1)!(k+1)}(x-x^*)^{k-1}$$\\\\
Which is equivalent to the expression of theorem 5 \quad $\Box$.
\\\\\\\\
\textbf{Remark:} With these results it is almost impossible for one to obtain a unstable shock curve in numerical simulation which may occur if one does not choose a initial point that satisfies these asymptotic relations. In the case when $h'(x)$ blows up at $x^*$ like $f(x)=\sqrt{x}$, $f(x)=xln(x)$ or even $f(x)=e^{\frac{1}{x}}$ in Burger's equation, one can still try the same techniques to obtain the asymptotic relations. For example, set $A>0$, $B>0$, $0<\alpha<1$, consider the Burger's equation with initial condition\\\\
$$
f(x) =
\begin{cases} 
A(-x)^{\alpha},  & \mbox{if }x\le0\\
-Bx^{\alpha}, & \mbox{if }x>0\\
\end{cases}
$$\\\\
It can be verified both by the mathematical techique or numerical simulation that $\frac{-\triangle x_l}{\triangle x_r}\sim s$ where s is the unique positive root of the function\\\\
$$s^{1+\alpha}+\frac{B(1+\alpha)}{A(1-\alpha)}s-\frac{1+\alpha}{1-\alpha}s^{\alpha}-\frac{B}{A}$$\\\\
In numerical simulations, it is unusal to apply such kind of initial condition. But if we must apply this method to these initial conditions, normally just setting $(0,x^*-\delta,x^*+\delta)$ for some small $\delta$ as initial point works very fine.
\newpage

\section{Numerical Method and Algorithm}
With all the content discussed in the previous sections, we are now well prepared to built an outline of the algorithm to solve (1.1.1). At least the math part is almost done, and it is now the time about the design of algorithm which must guarantee it obeys the mathematical facts and indeed solve what we expect. Of course as we will see there are still some details that one needs to use math to analyze but, theoretically, they are all solved and the rest are all technical works.
\subsection{Algorithm Logic and Stability Condition}
We are now going to discuss the case when $f(x)$ is a piecewise smooth function. The first step of the algorithm is certainly finding the shock points, determine which kind of shock points they are, and compute the needed derivatives of them. Then from each shock point, compute its break time. From the shock point which has the smallest break time, the algorithm starts to use some ODE solver to solve the shock curve. The following shows the logic of the basic algorithm design\\\\
$$
\begin{boxedminipage}{12.25cm}
\vspace{-15pt}
\begin{align*}
& \leftline{\textbf{Basic Characteristic Sweeping Algorithm}}\\\\
& \leftline{\text{1 \quad Input piecewise initial condition $f(x)$} and the desired time $T$ for solution}\\\\
& \leftline{\text{2 \quad Determine all potential shock / rarefaction / unsmooth points of $f(x)$}}\\
& \leftline{\text{\qquad and compute their side derivatives to the order used for asymptotic}}\\
& \leftline{\text{\qquad conditions}}\\\\
& \leftline{\text{3 \quad Start from the shock point that with the smallest break time, use ODE}}\\
& \leftline{\text{\qquad solver to solve its shock curve}}\\\\
& \leftline{\text{4 \quad Deal with the cases when shock point is cancled, shock curves merge,}}\\
& \leftline{\text{\qquad shock cross rarefaction, rarefaction cross rarefaction, and straight line}}\\
& \leftline{\text{\qquad connection if they occur}}\\\\
& \leftline{\text{5 \quad When all the shock curve before time $T$ are solved, start from initial}}\\
& \leftline{\text{\qquad condition, use characteristic to sweep out the solution at time $T$}}
\end{align*}
\vspace{-20pt}
\end{boxedminipage}
$$\\\\
We shall define the straight line connection in the next subsection. Step 1 and Step 2 are just trivial and not related to the content of this section. Now we must deal with all the cases start from step 3.\\\\
Before we start analyzing this process. let us discuss if this algorithm is stable. The stability of this method actually comes from the geometrical fact. As one may notice that given a fixed initial condition $f(x)$, all the shocks must come from all of the defined shock points, straightline connection, or their mergence (This can be mathematically proved but just some simple analysis of the previous definition). So at each time step when there is no case of condition change occur, we must have\\\\
$$t\le-\frac{1}{h'(x_l)} \qquad \text{and} \qquad t\le-\frac{1}{h'(x_r)}$$\\\\
for some $x_l$ and $x_r$ if $h'(x_l)<$ or $h'(x_r)<0$. Overall we mush have\\\\
$$1+h'(x)t\ge0$$\\\\
If we substitute $t$ with $t=-\frac{x_l-x_r}{h(x_l)-h(x_r)}$, then locally(so that there is no equal) when a shock forms, we have\\\\
$$h'(x_l)>\frac{h(x_l)-h(x_r)}{x_l-x_r} \qquad \text{and} \qquad h'(x_r)>\frac{h(x_l)-h(x_r)}{x_l-x_r} \eqno(4.3.1)$$\\\\
(4.3.1) is called \textbf{local stability condition}, and geometrically it can be treated as somewhat $x^*$ is a local "inflection point". In numerical simulations, one might obtain correct shock curve even without initial points that satisfy those asymptotic relation, but this initial points must satisfy (4.3.1).\\\\
Recall the shock equations we use, the most advantageous and important part of this method is: it does not only solve the shock curve coordinate $(\xi(t),t)$, but both the initial condition $x$-coordinates $x_l$ and $x_r$ whose function value $f(x_l)$ and $f(x_r)$ are indeed the solution value on the sides of the shock curve. This allow one to determine where to jump to the other side of the shock curve when characteristic intersect the shock curve during the process of sweeping.\\\\
The general sweeping process is described like this: given the solution time $T$, start from a small enough $x_L$ (Of course one can sweep from a large enough $x_R$ backward), compute the $x$ coordinate at time $T$ denote by $X_k=h(x_L+k\triangle x)T+(x_L+k\triangle x)$ for a small $\triangle x$, then the solution value for $X_k$ is $f(x_L+k\triangle x)$. If for some $k$ the $x_L+k\triangle x$ cross the first $x_l$ of some shock curve, then jump to its $x_r$ (If it is the shock adjacent to rarefaction then jump to the rarefaction point) and continue the sweeping. If $x_L+k\triangle x$ arrive some rarefaction point $x^*$, then compute $\frac{X-x^*}{T}$ and use (1.3.4) to sweep at time $T$ until $\frac{X-x^*}{T}$ cross out the region of rarefaction. Continue this whole process until the all $X$ coordinatex have coverd the $x$ span which our solution needs.\\\\
If for some $x_0$ and fixed $\triangle x$ such that $h'(x_0)$ is very large, then we can see\\\\
$$
\begin{cases}
X_0=h(x_0)T+x_0\\
X=h(x_0+\triangle x)T+x_0+\triangle x\\
\end{cases}
$$\\\\
Which gives $\triangle X=T\left[h(x_0+\triangle x)-h(x_0)\right]+\triangle x\sim(1+h'(x_0)T)\triangle x$ is also very large especially when T is large as well. In this case we can shrink the $\triangle x$ to $\triangle x=\frac{\triangle X}{1+h'(x_0)T}$ for a given expect $\triangle X$ and let $u(X,T)=f(x_0+\frac{\triangle X}{1+h'(x_0)T})$ be its approximate solution. Of course if one set the $\triangle x$ small enough this problem can be ignored.\\\\
\newpage

\subsection{Serveral Cases of Condition Change}
In this subsection we are going to numerically deal with the some cases that may occur in real algorithm implement.\\\\
\textbf{Case 1: Fake Shock Point}\\\\
As previously discussed a shock point does not necessarily form a shock. The reason is simple: some $x_l$ or $x_r$ may cross the shock point before its break time comes. In this case, what algorithm needs to do is just cancled this shock point and continue solving the shock curve. Note that for a fixed initial value $f(x)$, this cancellation is permanent.\\\\
Without loss of generality, we assume now there is a $x_r$ that is close to a shock point $x^*$ on its right side. At this specific time step, the algorithm detected a order change of $x_r$ and $x^*$. We can make it compute the derivatives of $x_r$ up to any order by using the equation (2.2.6) or (2.2.7). Then we can solve the following polynomial equation with respect to $\triangle t$ for a large enough $n$\\\\
$$x_r+\sum_{k=1}^{n}\frac{1}{k!}\frac{\mathrm{d}^kx_r}{\mathrm{d}t^k}\mid_{(x_l,x_r,t)}(\triangle t)^k=x^* \eqno(4.2.1)$$
Normally solving the first order or second order of this equation for a small $\triangle t$  is accurate enough, like $\triangle t=\frac{x^*-x_r}{\frac{\mathrm{d}x_r}{\mathrm{d}t}\mid_{(x_l,x_r,t)}}$. Then compare $t+\triangle t$ and the $t_b$ of $x^*$. If $t+\triangle t\le t_b$ the shock point $x^*$ is fake and will be cancled, otherwise locally we change the time step $\triangle t=t_b-t$ and get a new $x_r$ which is still on the left side of $x^*$. The next step is prepare for the shock mergence.\\\\
In fact there is a smarter way to deal with this case: at the beginning of one shock point, compute the the difference between its break time $t_{b_i}$ and the break time of the next shock point $t_{b_j}$. Then for a fix time step length $dt$, compute the next step solution for a small time step $\triangle t\le dt$ so that $\frac{t_{b_j}-t_{b_i}}{dt}$ is an integer, then use time step $dt$ back. To achieve this, just set $\triangle t=t_{b_j}-\lfloor\frac{t_{b_j}-t_{b_i}}{dt}\rfloor dt$. Here $\lfloor*\rfloor$ denote the sign of floor.\\\\\\
\textbf{Case 2: Shock Mergence}\\\\
Suppose there are 2 shock points $x_1^*<x_2^*$ and at one specific time step $x_{r_1}$ and $x_{l_2}$ are going to change their order. This time we are going to solve the following polynomial equation with respect to ${\triangle t}$ for a large enough n\\\\
$$x_{r_1}+\sum_{k=1}^{n}\frac{1}{k!}\frac{\mathrm{d}^kx_{r_1}}{\mathrm{d}t^k}\mid_{(x_{l_1},x_{r_1},t)}(\triangle t)^k=x_{l_2}+\sum_{k=1}^{n}\frac{1}{k!}\frac{\mathrm{d}^kx_{l_2}}{\mathrm{d}t^k}\mid_{(x_{l_2},x_{r_2},t)}(\triangle t)^k \eqno(4.2.2)$$\\\\
if n=1,then $\triangle t=\frac{x_{r_1}-x_{l_2}}{\frac{\mathrm{d}x_{l_2}}{\mathrm{d}t}\mid_{(x_{l_2},x_{r_2},t)}-\frac{\mathrm{d}x_{r_1}}{\mathrm{d}t}\mid_{(x_{l_1},x_{r_1},t)}}$. For such a small time step, solve the shock equation first, then at next step, $x_{1_r}$ and $x_{2_l}$ are merged so that cancled. The new shock equation will start with $(t+\triangle t, x_{l_1}\mid_{(t+\triangle t)},x_{r_2}\mid_{(t+\triangle t)})$.\\\\\\
\textbf{Case 3: Shock Cross Rarefaction}\\\\
In this case it is similar to the case 1, since the rarefaction point is still. The only differentce is the equation that are locally being solved is changed from the first kind to the second kind, or from the second kind to the first kind, depending on it is shock crossing into  or crossing out the rarefaction. To determine if a shock is crossing out a rarefaction, just check each time step if $\frac{\xi(t)-x^*}{t}$ still belongs to the region of rarefaction.\\\\\\
\textbf{Case 4: Rarefaction Cross Rarefaction}\\\\
It is trivial to see that if the two rarefaction points are not connect by a single straight line, then except for changing the equation from the second kind to the thrid kind, there is no difference between this case and the combination of Case 1 2 3. The function $f(x)$ can be anything between two adjacent rarefaction points. Just let the algorithm detemine which case it is. Here is a simple property for such a shock appear in the intersection of rarefaction.\\\\
\textbf{Property 2:} \emph{A shock created by two rarefaction intersection will not cross out the region of these rarefactions unless it merges with other shocks that cross in the rarefaction region.}\\\\
This property is obvious as one can see from the third kind shock equation, the speed of such a shock must between the speed of the two boundary rays of rarefaction region.\\\\\\
\textbf{Case 5: Straight Line Connection}\\\\
If understand fully about the defined kinds of shock point so far, one should have noticed that the the discussion of one case is missing\\\\
\textbf{Definition 5:} \emph{A point $x^*$ is called \textbf{straight line connection} if the following conditions is satisfied}\\\\
1 \quad \emph{f(x) is continuous at $x^*$}\\\\
2 \quad \emph{There exists $\delta_0>0$ such that for any $\delta\in(0,\delta_0)$, at least one of the follwoing \emph{(1) (2)} is satisfied}\\\\
\emph{\emph{(1)} $h(x)=k(x+c)$ for any $x\in\left(x^*-\delta,x^*\right]$ and some k such that $k<0$ and $k\le h'(x^{*+})$}\\\\
\emph{\emph{(2)} $h(x)=k(x+c)$ for any $x\in\left[x^*,x^*+\delta\right)$ and some k such that $k<0$ and $k\le h'(x^{*-})$}\\\\\\
This however does not mean the method cannot handle it. The reason to discuss this independently is the fact that in the case when $h(x)$ is straight line, its characteristic all intersect to a single point, and locally the solution of $u(x,t)$ can be explicitly written down. Assume now $h(x)=k(x+c)$ at interval $\left[a,b\right]$, according to (2.1.1) we have\\\\
$$
\begin{cases}
k(x_l+c)t=x-x_l\\
k(x_r+c)t=x-x_r\\
\end{cases}
$$\\\\
Which can be solve as $t=-\frac{1}{k}$ irrelevant to $x_l$ and $x_r$. Moreover according to (1.1.8), assume $u(x,t)=z=f(x_0)$ for some $x_0\in\left[a,b\right]$, then\\\\
$$f(x_0)=f(x-G'(f(x_0))t)=f(x-h(x_0)t)=f(x-k(x_0+c)t)$$\\\\
According to the definition of straight line connection, $f$ must be monotonic at interval $\left[a,b\right]$. This means $x_0=x-k(x_0+c)t$ which is\\\\
$$x_0=\frac{x-kct}{1+kt}$$\\\\
Thus, the solution of $u(x,t)$ in the region $\{(x,t)\mid 0\le t<\frac{1}{k}, x\in\left[k(a+c)t+a,k(b+c)t+b\right]\}$ is\\\\
$$u(x,t)=f(\frac{x-kct}{1+kt}) \eqno(4.2.3)$$\\\\
(4.2.3) can be sometimes used to improve accuracy of solving (2.2.6) when one side is a straight line. \textbf{If $x_1^*<x_2^*$ are two adjacent straight line connection and the tangent of its straight line is $k<0$, then if when $t<-\frac{1}{k}$ there is no any points $x_l$ and $x_r$ of other shock points coming in the interval $(x_1^*,x_2^*)$, then cancle these two straight line connections and immediately start a shock from the left side of $x_1^*$ and right side of $x_2^*$.}\\\\\\
\textbf{Case 6: Local Differential Singularity at Discontinous Point}\\\\
In the case 2 or case 3, if at time step t the shock equation changes or shock equation arrive to a point which is a singular point of equation $(2.2.6)$ or $(2.2.7)$, then certainly we need local analysis of the equations which is different from what we have done, since it is not at the initial time. However there are not much difference: when it is going to change to case 3, the previous theorem 5 is still applicable; when it is going to start a new shock whose at least one side's derivative is singular, we can use the techniques we used with a little adjustment. Recall we still have\\\\
$$\frac{\mathrm{d}x_l}{\mathrm{d}x_r}=A(x_l,x_r)\cdot B(x_r,x_r)$$\\\\
where\\\\
$$A(x_l,x_r)=\frac{G(f(x_r))-G(f(x_l))-G'(f(x_l))(f(x_r)-f(x_l))}{G(f(x_r))-G(f(x_l))-G'(f(x_r))(f(x_r)-f(x_l))}$$\\\\
$$B(x_l,x_r)=\frac{1+G''(f(x_r))f'(x_r)t}{1+G''(f(x_l))f'(x_l)t}$$\\\\
For the analysis of $B(x_l,x_r)$ is the same as that previous used. Now if the new shock is continuous at $f(x_r)$ and $f(x_l)$ just like the case in which we previously discuss, then we still have $A(x_l,x_r)\to-1$ and theorem 1 2 3 4 are all applicable; If $f(x_l)\ne f(x_r)$, then $A(x_l,x_r)$ is not singular and can be computed exactly. Substitute this value to the previous analysis process of theorem 1 2 3 4, then we can get the asymptotic relation for this case. For example, suppose the new shock has $A(x_l,x_r)=-2$, then (3.1.2) becomes\\\\
$$-2s=\frac{(h_+^{(k)}-h_-^{(k)}s^k)-kh_+^{(k)}(1-s)}{(h_+^{(k)}-h_-^{(k)}s^k)-kh_-^{(k)}(1-s)s^{k-1}}$$\\\\
Therefore this case can be solved just like how we solve the previous case. And there are certainly new formulas of theorem 2 3 4 which contain the term $A(x_l,x_r)$. We denote $A=A(x_l,x_r)$ , then for theorem 2, the polynomial equation becomes\\\\
$$s^{k+1}+\frac{1-Ak}{A(k-1)}s^k+\frac{A-k}{A(k-1)}ps+\frac{p}{A}=0$$\\\\
For theorem 3, the asymptotic relation becomes\\\\
$$\triangle x_l\sim\frac{(k_r-1)h_+^{(k_r)}}{(1-Ak_r)k_r!(h_+'-h_-')}\triangle x_r^{k_r}$$\\
$$\triangle x_r\sim\frac{(k_l-1)h_+^{(k_l)}}{(1-Ak_l)k_l!(h_-'-h_+')}\triangle x_l^{k_l}$$\\\\
Similarly for theorem 4\\\\
$$\triangle x_l\sim-\left|\left[\frac{(1-k_r)k_l!h_+^{(k_r)}}{(1-Ak_r)k_r!h_-^{(k_l)}}\right]^{\frac{1}{k_l}}\right|\triangle x_r^{\frac{k_r}{k_l}}$$\\
$$\triangle x_r\sim\left[(-1)^{k_l}\frac{(1-k_l)k_r!h_-^{(k_l)}}{(1-Ak_l)k_l!h_+^{(k_r)}}\right]^{\frac{1}{k_r}}(-\triangle x_l)^{\frac{k_l}{k_r}}$$\\\\\\
\textbf{Case 7: Combination of All Cases Above}\\\\
If at a specific time step more than one of the cases above occur, than just compute each of their $\triangle t$ and choose a smallest one as the next time step length. Moreover, for case 3 when the start derivatives is singular, like some of the cases of theorem 5, we can locally flip the equation and solve $\frac{\mathrm{d}t}{\mathrm{d}x}$ first for some small $x$ steps then flip it back.\newpage

\subsection{Numerical Simulations and Comparision with Traditional Methods}
Now we are going to compare CS methods with those traditional methods.\\\\\\
\textbf{Operation time}\\\\
The main advantage of this CS method compare to those traditional methods is the operation time. Traditional methods require $O(mn)$ time where $m$ is the amout of $x$ coordinate partition and $n$ is the amout of $t$ coordinate partition. CS method only needs $O(kn)$ time where $k$ is the total amount of shock / rarefaction points and straight line connections. In real simulations the operation time are expected to be less than $O(kn)$ since there may be some cases when some of these points are cancled.\\\\\\
\textbf{Storage}\\\\
Moreover because in CS method the solution can be swept out if we have those $x_l$ and $x_r$, if one needs to store the dynamic solution, traditional methods need $O(mn)$ storage and CS method only needs $O(kn)$ storage as well. Since normally $k<<m$, this method should work much better in most of the cases. \\\\\\
\textbf{Stability}\\\\
In most of the cases I have simulated, if one use the initial points that obey those asymptotic relation then the shock curve is very stable and correct. Compared to the traditional methods which normally have CFL condition, CS method does not have restriction on time step length.\\\\\\
\textbf{Accuracy}\\\\
The accuracy / error depends on the chosen initial points and the accuracy of the ODE solver. The asymptotic relation that used to choose initial points can be arbitrarily generalized to high order by using the same techniques (See the reference). And the ODE solver such as Runge-Kutta method, theoretically can go to arbitrarily high order as well. So the accuracy of CS method is guaranteed as we will see at the following pictures of numerical simulations.\newpage
\noindent\textbf{Note: All of the following examples are simulated with 5-th order WENO scheme and CS method with 4-th order Runge-Kutta ODE solver. Because the sweeping process is exact, one only needs to test how well the CS method predicts the shock curve. The background is the 3D graph generated by WENO and the red curve is generated by CS. CPU time only records for computation process.}\\\\
\textbf{Example 1: First Kind Shock Point}\\\\
$$G(u)=\frac{1}{2}u^2 \qquad \qquad f(x)=\begin{cases} 
x+1.5,  & \mbox{if }x<0\\
x^2-2x, & \mbox{if }x>0\\
\end{cases}$$
\par
\centerline{\includegraphics[width=20cm,height=13cm]{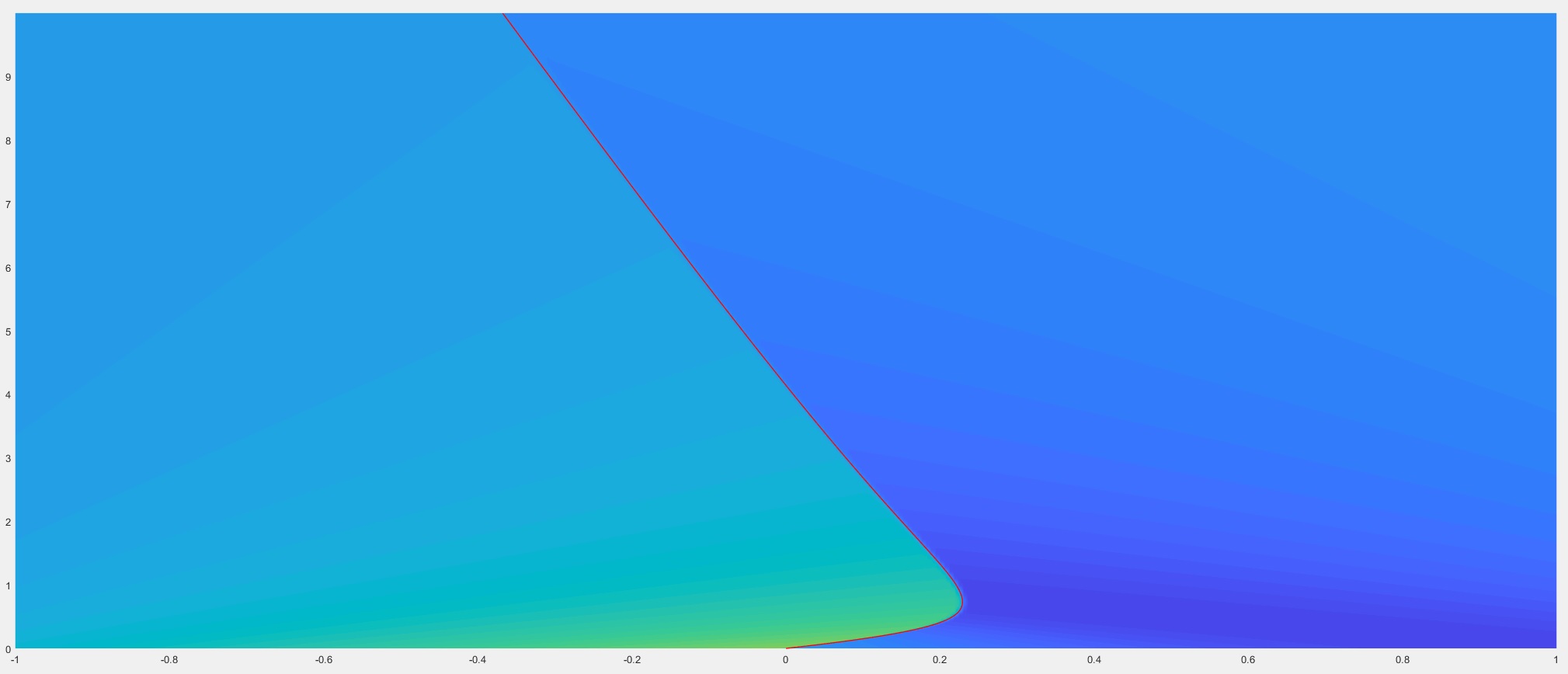}}
\centerline{Example 1 \qquad dx=0.006, dt=0.0012, T=10}
\centerline{CPU time: WENO 1.354s; CS 0.351s}
\par
\centerline{\includegraphics[width=24cm,height=16cm]{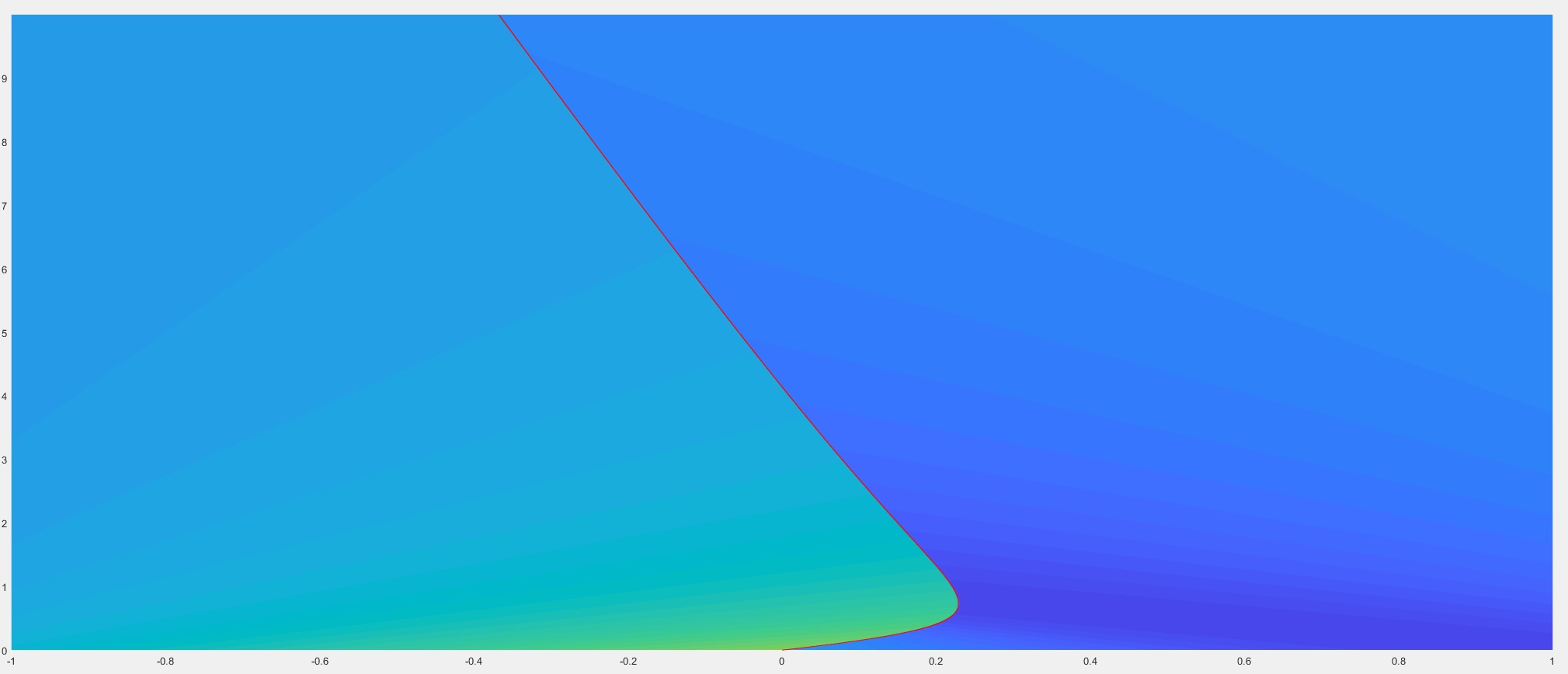}}
\centerline{Example 1 \qquad dx=0.0012, dt=0.00024, T=10}
\centerline{CPU time: WENO 95.79s; CS 1.648s}\newpage

\textbf{Example 2: Second Kind Shock Point}
$$G(u)=\frac{1}{2}u^2 \qquad \qquad f(x)=\begin{cases} 
1-e^x,  & \mbox{if }x<0\\
x^2-2x, & \mbox{if }x>0\\
\end{cases}$$
\par
\centerline{\includegraphics[width=24cm,height=16cm]{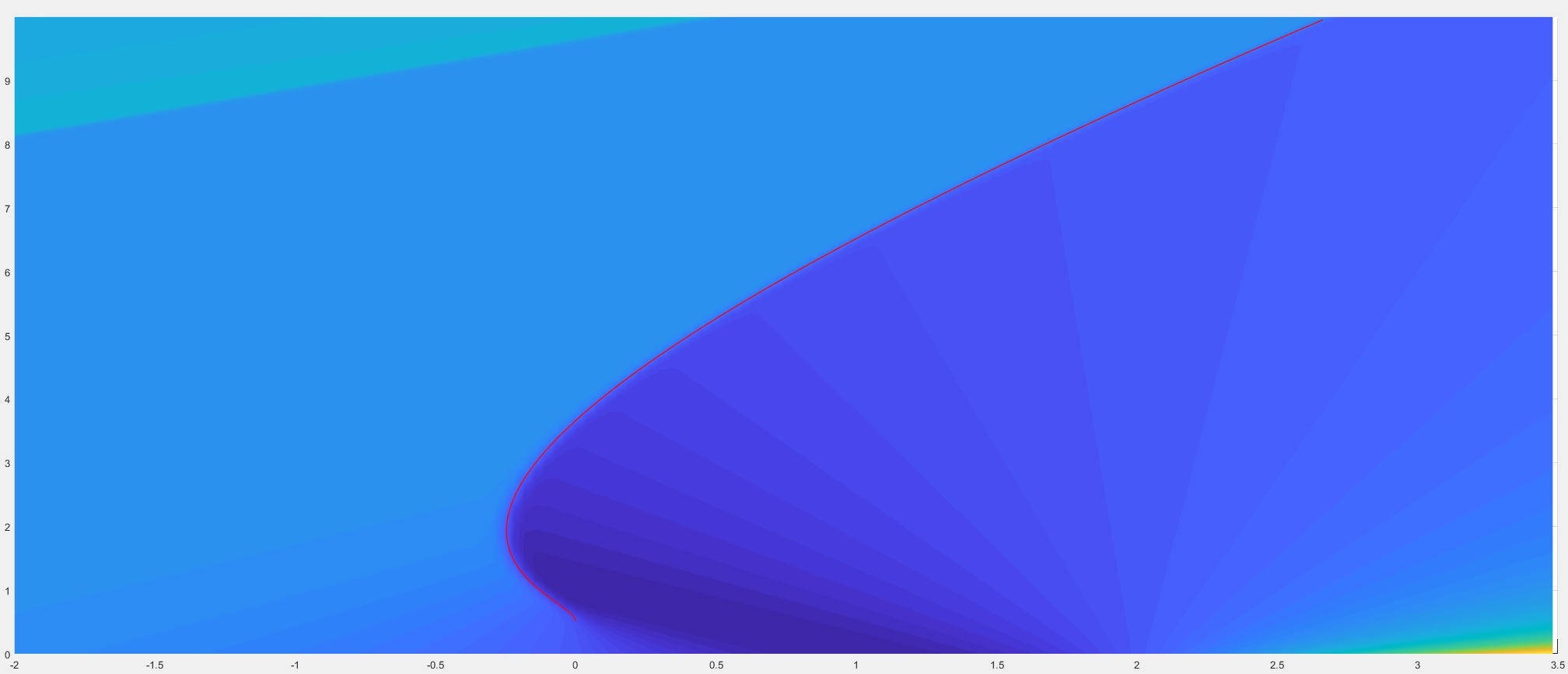}}
\centerline{Example 2 \qquad dx=0.037, dt=0.0067, T=10}
\centerline{CPU time: WENO 0.217s; CS 0.0295s}
\par
\centerline{\includegraphics[width=24cm,height=16cm]{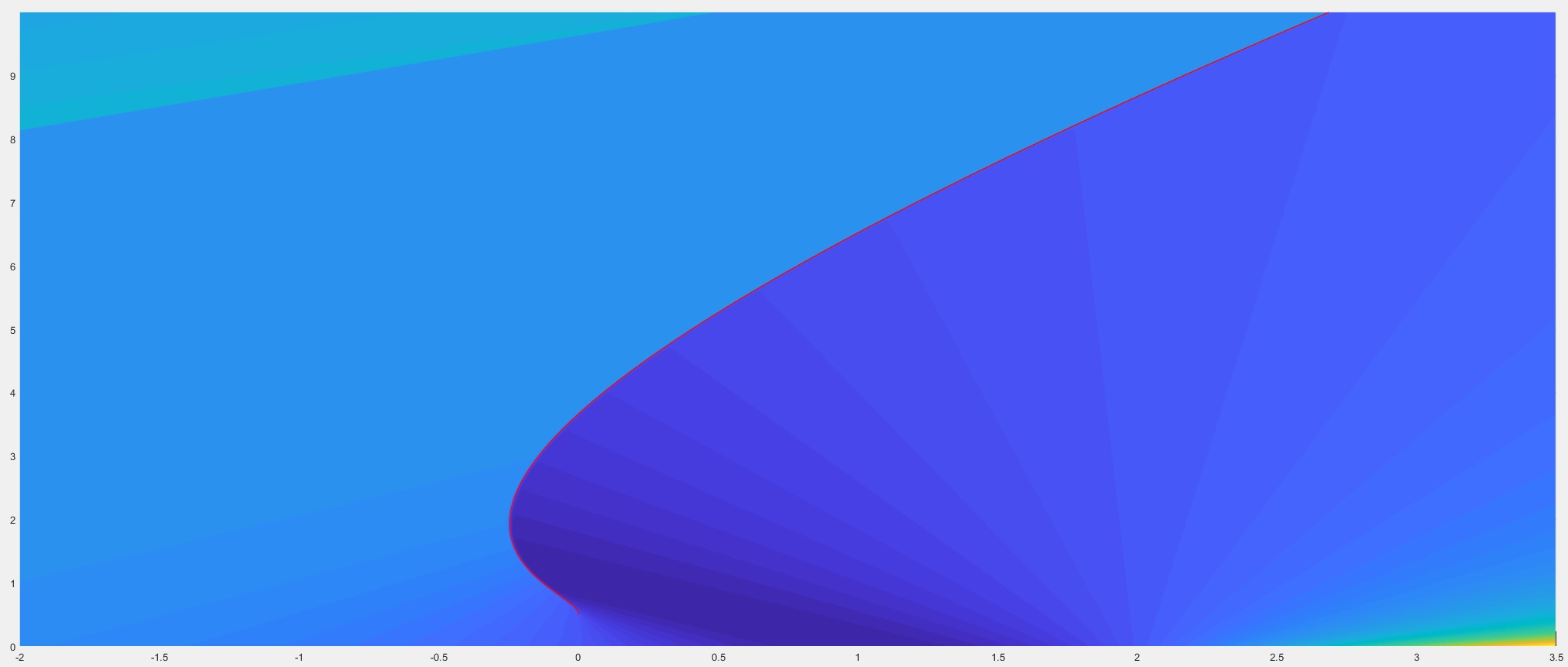}}
\centerline{Example 2 \qquad dx=0.037, dt=0.0067, T=10}
\centerline{CPU time: WENO 4.99s; CS 0.0375s}\newpage

\textbf{Example 3: Third Kind Shock Point}\\\\
$$G(u)=\frac{1}{2}u^2\qquad \qquad f(x)=\frac{2x}{(1+x^2)^2}$$\\\\
Use the exact solution of Burger's equation we can actually in this case get the analytic solution of the shock curve: $t=\frac{1}{2}x^2$ and $t=-\frac{1}{2}x^2$.\\\\
\par
\centerline{\includegraphics[width=24cm,height=14cm]{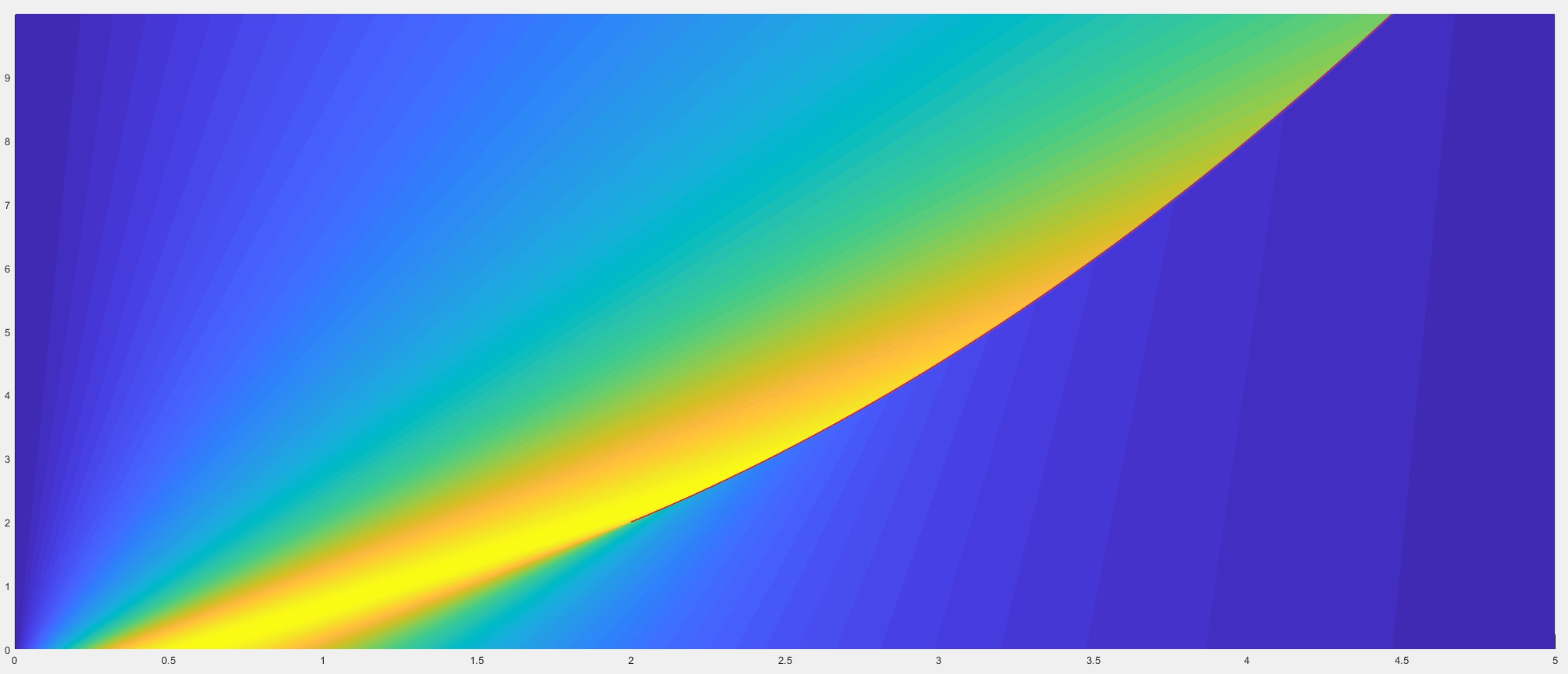}}
\centerline{Example 3 \qquad dx=dt=0.005, T=10}
\centerline{CPU time: WENO 0.356s; CS 0.085s}\newpage

\textbf{Example 4: Fourth Kind Shock Point}
$$G(u)=\frac{1}{12}u^4 \qquad \qquad f(x)=\begin{cases} 
\frac{-x^2-2x-1}{2},  & \mbox{if }x\le0\\
x+1, & \mbox{if }x>0\\
\end{cases}$$
\par
\centerline{\includegraphics[width=24cm,height=16cm]{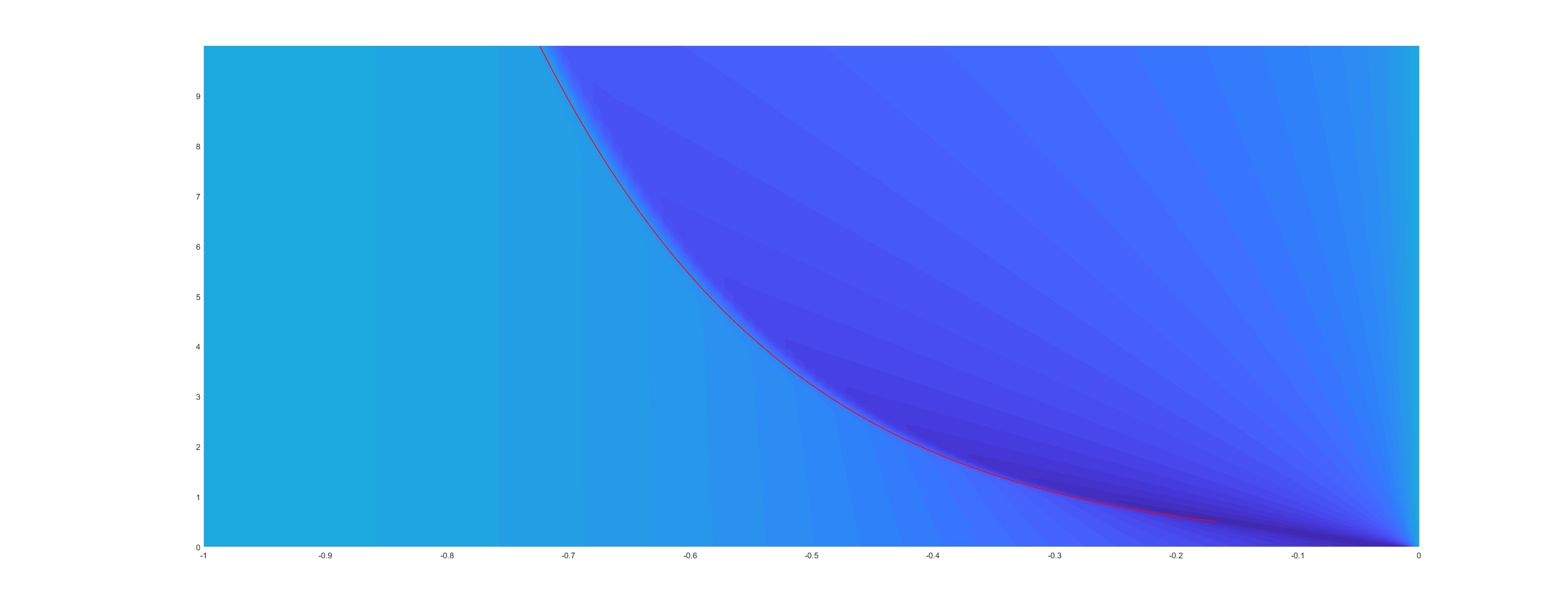}}
\centerline{Example 4 \qquad dx=0.0833, dt=0.01, T=10}
\centerline{CPU time: WENO 3.1s; CS 0.052s}
\par
\centerline{\includegraphics[width=24cm,height=16cm]{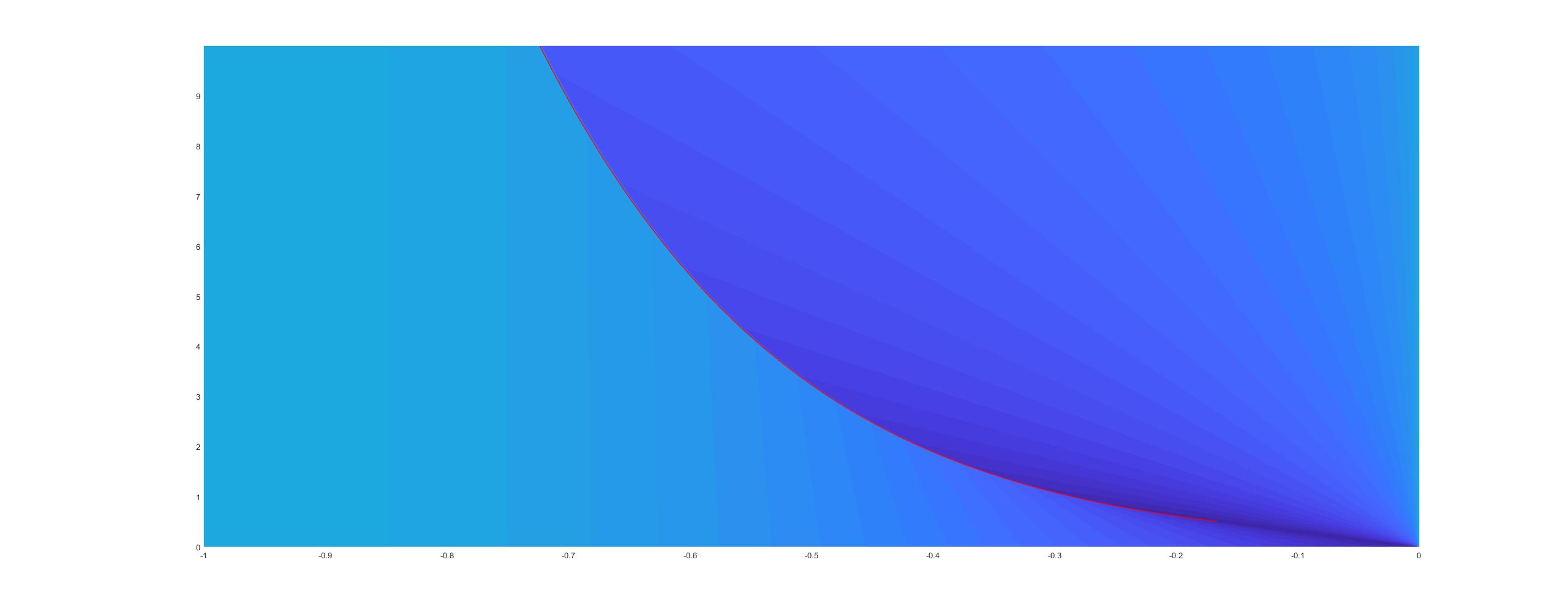}}
\centerline{Example 4 \qquad dx=0.00833, dt=0.001, T=10}
\centerline{CPU time: WENO 103.2s; CS 0.392s}\newpage

\textbf{Example 5: Shock Mergence}
$$G(u)=\frac{1}{2}u^2\qquad \qquad f(x)=e^{\frac{-x^4+5x^2}{10}}$$\\\\
\par
\centerline{\includegraphics[width=24cm,height=15cm]{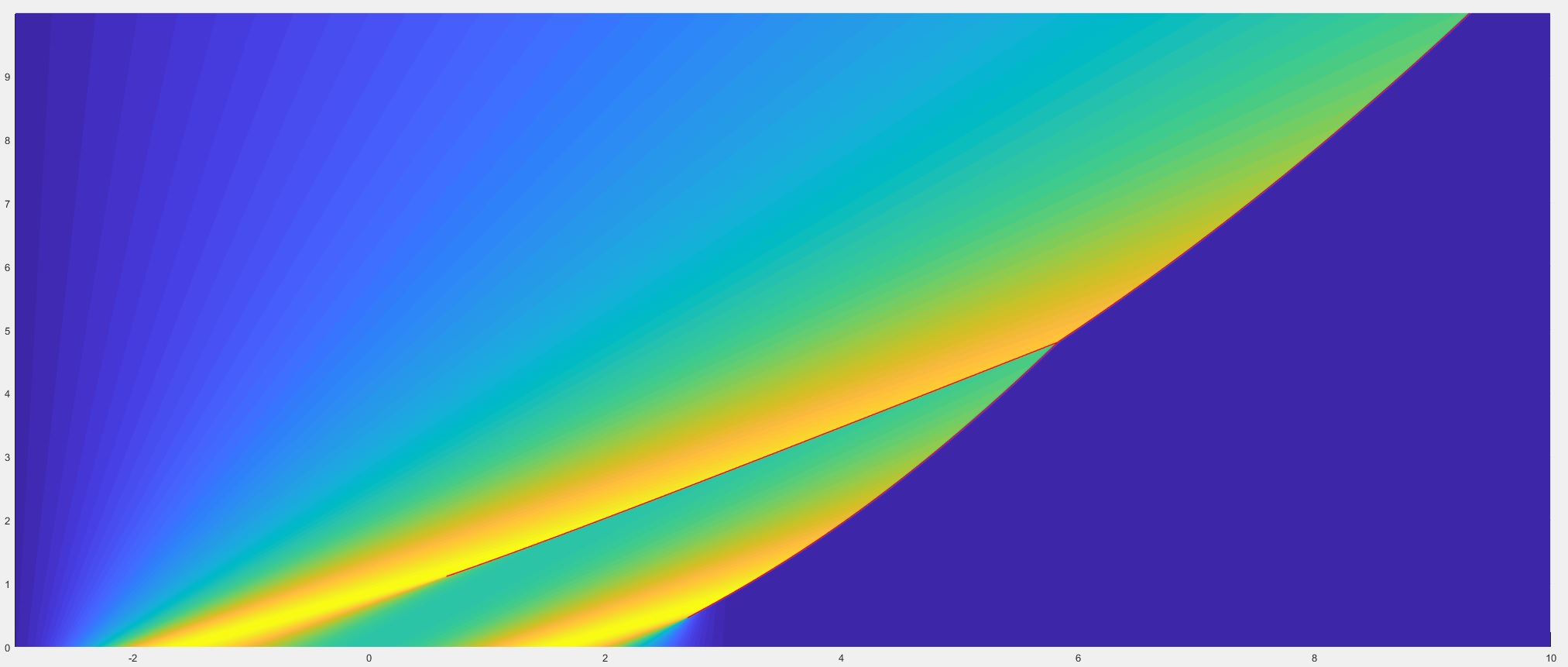}}
\centerline{Example 5 \qquad dx=0.026, dt=0.013, T=10}
\centerline{CPU time: WENO 0.288s; CS 0.064s}
\par
\centerline{\includegraphics[width=24cm,height=16cm]{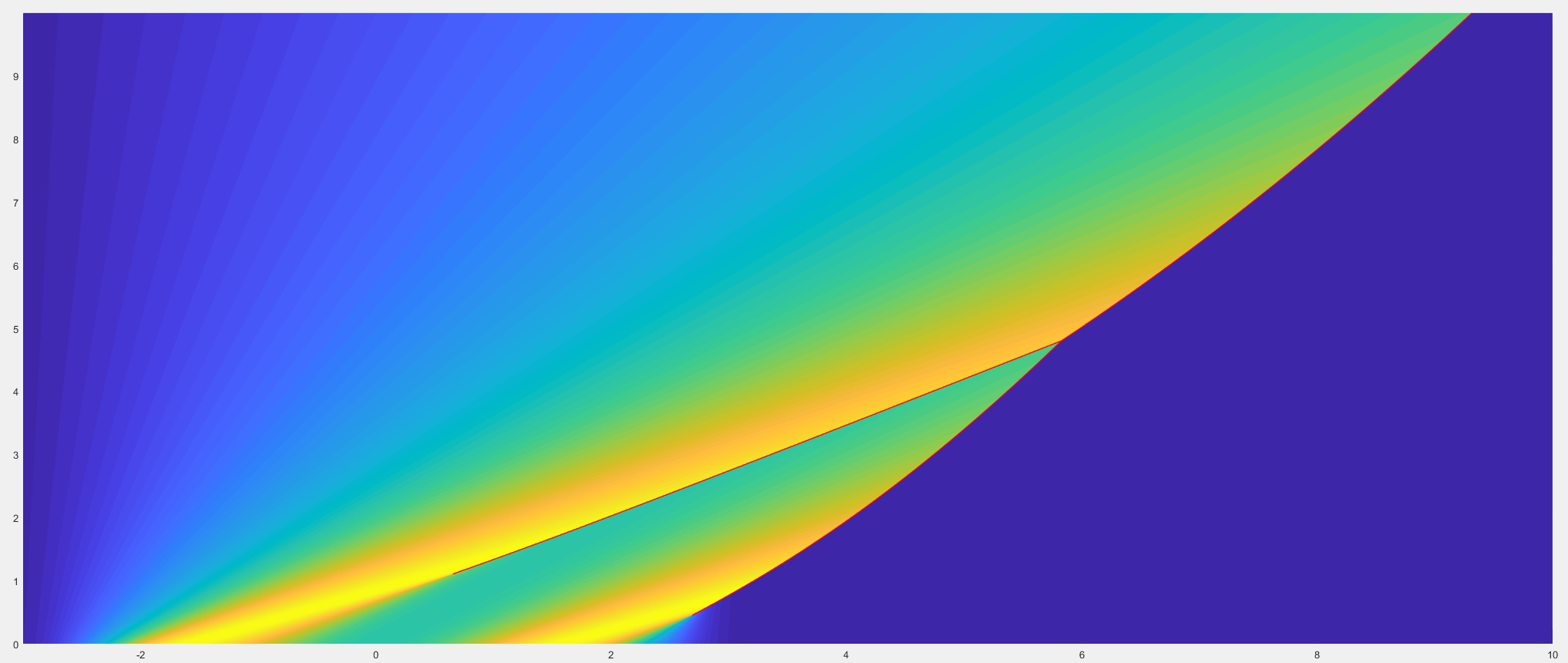}}
\centerline{Example 5 \qquad dx=0.0026, dt=0.0013, T=10}
\centerline{CPU time: WENO 57s; CS 0.45s}\newpage

\noindent\textbf{Remark:} From the numerical results we can see, CS method is very accurate compared with traditional methods even it runs at a low time step length. In fact this is expected since once it is stable, it is solving the shock curve equation exactly. Moreover one can also numerically verify those asymptotic relations and see they are indeed true. And what may make one surprised is we sometimes get a correct shock curve even we don't use initial points that obey the asymptotic relations.\\\\
CS method now only works in 1D. It seems in real applications no one would like to solve 1D problem. But in fact among all kinds of the numerical methods on numerically solving PDE, all of the advantages of CS method along with its restriction on 1D can be used as the tool of error analysis of the methods which use traditional time step iteration for all $x$ span, since if one comes up with a new method, normally it will be applied this to 1D problem first to see how well it works.\\\\
The biggest disadvantage of CS method is also obvious: it is hard to be coded. It involves so many cases management and logic judgement, and needs a specific data structure. So far I just mostly simulate cases locally, and I have not even started to code this method. However, I am not a prefessional programmer but I really hope some can code this out.\\\\\\\\\\\\

\end{document}